\documentclass[12pt,a4paper]{article}
\usepackage{amsmath}
\usepackage{amsthm}
\usepackage{graphicx}
\allowdisplaybreaks[1]
\newcommand{\os}{\underline}

\newcommand{\pr}{\rightarrow}

\newcommand{\ba}{\begin{array}}
\newcommand{\ea}{\end{array}}

\newcommand{\varp}{\varphi}
\newcommand{\eps}{\varepsilon}

\newcommand{\il}{\int\limits}
\newenvironment{inspring}[1]%
{\begin{list}{}{\setlength{\rightmargin}{0cm}
                \setlength{\listparindent}{0cm}
                \settowidth{\labelwidth}{\mbox{#1}}
                \setlength{\leftmargin}{1.1\labelwidth}
                \setlength{\labelsep}{.1\labelwidth}}}%
{\end{list}}
\newcommand{\ITEM}[1]{\item[#1\hfill]}
\newcommand{\bi}[1]{\begin{inspring}{#1}}
\newcommand{\ei}{\end{inspring}}

\newcommand{\dlim}{\displaystyle \lim}

\newcommand{\dsum}{\displaystyle \sum}

\newcommand{\beq}{\begin{equation}}
\newcommand{\eq}{\end{equation}}
\catcode`\@=11

\font\tenmsa=msam10 \font\sevenmsa=msam7 \font\fivemsa=msam5
\font\tenmsb=msbm10 \font\sevenmsb=msbm7 \font\fivemsb=msbm5
\newfam\msafam
\newfam\msbfam
\textfont\msafam=\tenmsa  \scriptfont\msafam=\sevenmsa
  \scriptscriptfont\msafam=\fivemsa
\textfont\msbfam=\tenmsb  \scriptfont\msbfam=\sevenmsb
  \scriptscriptfont\msbfam=\fivemsb

\def\Bbb{\ifmmode\let\next\Bbb@\else
 \def\next{\errmessage{Use \string\Bbb\space only in math mode}}\fi\next}
\def\Bbb@#1{{\Bbb@@{#1}}}
\def\Bbb@@#1{\fam\msbfam#1}
\newcommand{\dR}{{\Bbb R}}

\newtheorem{thm}{Theorem}
\newtheorem{lem}[thm]{Lemma}

\newtheorem{prop}[thm]{Proposition}

\title{Bounding Taylor approximation errors for the exponential function in the presence of a power weight function}
\author{A.J.E.M.\ Janssen \\
Eindhoven University of Technology \\
Department of Mathematics and Computer Science}
\date{}

\begin{document}
\maketitle
\mbox{} \\ \\ \\ \\ \\
\noindent
{\bf Abstract.} \\
Motivated by the needs in the theory of large deviations and in the theory of Lundberg's equation with heavy-tailed distribution functions, we study for $n=0,1,...$ the maximization of
\begin{eqnarray}
& \mbox{} &
S:~\Bigl(1-e^{-s}\Bigl(1+\frac{s^1}{1!}+...+\frac{s^n}{n!}\Bigr)\Bigr)/s^{\delta} = E_{n,\delta}(s)~{\rm over}~s\geq0 \nonumber \\[2mm]
& & \hspace*{1cm}{\rm with}~\delta\in(0,n+1)~, \nonumber \\[4mm]
& & U:~({-}1)^{n+1}\Bigl(e^{-u}-\Bigl(1-\frac{u^1}{1!}+...+({-}1)^n\,\frac{u^n}{n!} \Bigr)\Bigr)/u^{\delta}=G_{n,\delta}(u)~{\rm over}~u\geq0 \nonumber \\[2mm]
& & \hspace*{1cm}{\rm with}~\delta\in(n,n+1)~. \nonumber
\end{eqnarray}
We show that $E_{n,\delta}(s)$ and $G_{n,\delta}(u)$ have a unique maximizer $s=s_n(\delta)>0$ and $u=u_n(\delta)>0$ that decrease strictly from $+\infty$ at $\delta=0$ and $\delta=n$, respectively, to 0 at $\delta=n+1$. We use Taylor's formula for truncated series with remainder in integral form to develop a criterion to decide whether a particular smooth function $S(\delta)$, $\delta\in(0,n+1)$, or $U(\delta)$, $\delta\in(n,n+1)$, respectively, is a lower/upper bound for $s_n(\delta)$ and $u_n(\delta)$, respectively. This criterion allows us to find lower and upper bounds for $s_n$ and $u_n$ that are reasonably tight and simple at the same time. As a result, the maximum location of $E_{n,\delta}(s)$ and $G_{n,\delta}(u)$, as well as their maximum values $ME_{n,\delta}$ and $MG_{n,\delta}$, are accurately estimated. Furthermore, as a consequence of the identities $\frac{d}{d\delta}\,[{\rm ln}\,ME_{n,\delta}] ={-}{\rm ln}\,s_n(\delta)$ and $\frac{d}{d\delta}\,[{\rm ln}\,MG_{n,\delta}]={-}{\rm ln}\,u_n(\delta)$, we show that $ME_{n,\delta}$ and $MG_{n,\delta}$ are log-convex functions of $\delta\in(0,n+1)$ and $\delta\in(n+1,n)$, respectively, with limiting values 1 ($\delta\downarrow0$) and $1/(n+1)!$ ($\delta\uparrow n+1$) for $E$, and $1/n!\,(\delta\downarrow n)$ and $1/(n+1)!\,(\delta\uparrow n+1)$ for $G$. The minimal values $\hat{E}_n$ and $\hat{G}_n$ of $ME_{n,\delta}$ and $MG_{n,\delta}$, respectively, as a function of $\delta$, as well as the minimum locations $\delta_{n,E}$ and $\delta_{n,G}$ are determined in closed form.
\mbox{} \\ \\
\section{Introduction} \label{sec1}
\mbox{} \\[-9mm]

In Nagaev's survey paper \cite{ref1} on large deviations of independent random variables, a basic technical step is to bound a quantity
\beq \label{e1}
\il_X\,\Bigl(e^{tx}-1-tx-...-\frac{1}{n!}\,(tx)^n\Bigr)\,dF(x)
\eq
in terms of a moment
\beq \label{e2}
\il_{-\infty}^{\infty}\,|x|^{\delta}\,dF(x)~.
\eq
Here $t>0$, $n=0,1,...\,$, $F(x)$, $x\in\dR$, is a probability distribution function on $\dR$, $X$ is a subinterval of $\dR$ and $\delta>0$. This is achieved by writing the quantity in (\ref{e1}) as
\beq \label{e3}
\il_X\,\frac{e^{tx}-1-tx-...-\dfrac{1}{n!}\,(tx)^n}{|tx|^{\delta}}\,|tx|^{\delta}\, dF(x)~,
\eq
and to use a bound for the function
\beq \label{e4}
\frac{e^{tx}-1-tx-...-\dfrac{1}{n!}\,(tx)^n}{|tx|^{\delta}}~.
\eq
When $s=tx>0$, one writes
\beq \label{e5}
\frac{e^s-1-s-...-\dfrac{1}{n!}\,s^n}{s^{\delta}}=e^s\, \frac{1-\Bigl(1+s+...+\dfrac{1}{n!}\,s^n\Bigr)\,e^{-s}}{s^{\delta}}~,
\eq
and then the interest is in bounding
\beq \label{e6}
E_{n,\delta}(s)=\frac{1-\Bigl(1+s+...+\dfrac{1}{n!}\,s^n\Bigr)\,e^{-s}}
{s^{\delta}}~,~~~~~~s\geq0~.
\eq
For instance, in \cite{ref1}, (1.12--13) of the proof of Lemma~1.4, this is considered for $\delta=2$, $n=1$, and there the bound
\beq \label{e7}
\frac{1-(1+s)\,e^{-s}}{s^2}\leq1/2~,~~~~~~s\geq0~,
\eq
is used. Note that $E_{n,\delta}(s)$ is bounded in $s\geq0$ if and only if $\delta \in[0,n+1]$.

When $s=tx<0$, the interest is in bounding
\beq \label{e8}
G_{n,\delta}(u)=({-}1)^{n+1}\,\frac{e^{-u}-1+u-...-({-}1)^n\,\dfrac{u^n}{n!}}
{u^{\delta}}~,~~~~~~u={-}s\geq0~.
\eq
The factor $({-}1)^{n+1}$ in (\ref{e8}) has been included to achieve that $G_{n,\delta}(u)\geq0$, $u\geq0$. Note that $G_{n,\delta}(u)$ is bounded in $u\geq0$ if and only if $\delta\in[n,n+1]$.

Bounding $G_{n,\delta}(u)$ is also of interest when one considers Lundberg's equation, see, for instance, \cite{ref2}, for heavy-tailed distributions. Here one has a probability distribution function $F(x)$, $x\in\dR$, such that for some $\beta>1$
\beq \label{e9}
\il_{-\infty}^{\infty}\,|x|^{\beta}\,dF(x)<\infty~,
\eq
while $\il_{-\infty}^{\infty}\,\exp(\gamma x)\,dF(x)=\infty$ for all $\gamma>0$ (heavy-tailedness) and the expectation value $\mu=\il_{-\infty}^{\infty}\,x\,dF(x)$ is negative. One considers for $y\geq0$ and $\gamma\geq0$ the functional
\beq \label{e10}
L(y\,;\,\gamma):=\il_{-\infty}^y\,e^{\gamma x}\,dF(x)+e^{\gamma y}\,\il_y^{\infty}\,dF(x)~.
\eq
Due to $\mu<0$, the equation $L(y\,;\,\gamma)=1$ has for any $y\geq0$ (along with the trivial solution $\gamma=0$) a unique solution $\gamma=\gamma(y)>0$. We have $\gamma(y)\downarrow0$ as $y\pr\infty$, due to heavy-tailedness.

The equation $L(y\,;\,\gamma)=1$ with $y>0$ is known as Lundberg's equation, and below we sketch an approach to approximate its solution $\gamma(y)$ as $y\pr\infty$. One first rewrites $L(y\,;\,\gamma)$ for $\gamma>0$ as
\beq \label{e11}
L(y\,;\,\gamma)=\il_{-\infty}^0\,e^{\gamma x}\,dF(x)+(1-F(0))+\il_0^y\,\gamma\,e^{\gamma x}(1-F(x))\,dx~.
\eq
Next, one uses the inequality
\beq \label{e12}
0<\frac{e^{-u}-(1-u)}{u^{\delta}}\leq\frac{1}{\delta}~,~~~~~~u>0~.
\eq
with $\delta=\min\,\{\beta,2\}\in(1,2]$ to find
\begin{eqnarray} \label{e13}
\il_{-\infty}^0\,e^{\gamma x}\,dF(x) & = & \il_{-\infty}^0\,(1+\gamma x)\,dF(x)+\il_{-\infty}^0\,(e^{\gamma x}-1-\gamma x)\,dF(x) \nonumber \\[3mm]
& = & F(0)+\gamma\,\il_{-\infty}^0\,x\,dF(x)+\eps(\gamma)~,
\end{eqnarray}
where
\beq \label{e14}
0<\eps(\gamma)=\il_{-\infty}^0\,\frac{e^{\gamma x}-1-\gamma x}{|\gamma x|^{\delta}}\, |\gamma x|^{\delta}\,dF(x)\leq\frac{\gamma^{\delta}}{\delta}\,\il_{-\infty}^0\,|x|^{\delta}\,dF(x) =O(\gamma^{\delta})~.
\eq
Then it follows from (\ref{e11}) and $L(y\,;\,\gamma)=1$ that
\beq \label{e15}
1=F(0)+\gamma\,\il_{-\infty}^0\,x\,dF(x)+O(\gamma^{\delta})+(1-F(0))+\il_0^y\, \gamma\,e^{\gamma x}(1-F(x))\,dx~,
\eq
i.e., by simplifying and dividing through by $\gamma$, that
\beq \label{e16}
\il_0^y\,e^{\gamma x}(1-F(x))\,dx={-}\,\il_{-\infty}^0\,x\,dF(x)+O(\gamma^{\delta-1})~.
\eq
Finally, using that
\beq \label{e17}
\mu=\il_{-\infty}^0\,x\,dF(x)+\il_0^{\infty}\,(1-F(x))\,dx~,
\eq
one gets
\beq \label{e18}
\il_0^y\,e^{\gamma x}(1-F(x))\,dx={-}\mu+\il_0^{\infty}\,(1-F(x))\,dx+O(\gamma^{\delta-1})~.
\eq
A strategy to proceed from this point onwards is to ignore the $O(\gamma^{\delta-1})$ at the right-hand side of (\ref{e18}) (remembering that $\gamma(y)\pr0$ as $y\pr\infty)$, and to concentrate on the simplified equation under certain decay assumptions on $1-F(x)$, $x\pr\infty$.

The subject matter of the present article is therefore to bound for $n=0,1,...$ and $\delta\in[0,n+1]$ the quantity $E_{n,\delta}(s)$ over $s\geq0$, and to bound for $n=0,1,...$ and $\delta\in[n,n+1]$ the quantity $G_{n,\delta}(u)$ in (\ref{e8}) over $u\geq0$. We summarize the main results in the next section.

\section{Main results and overview} \label{sec2}
\mbox{} \\[-9mm]

We consider for $n=0,1,...$ and $\delta\in[0,n+1]$ the quantity
\beq \label{e19}
E=E_{n,\delta}(s)=\frac{1-\Bigl(1+\dfrac{s^1}{1!}+...+\dfrac{s^n}{n!}\Bigr)\,e^{-s}} {s^{\delta}}~,~~~~~~s\geq0~,
\eq
with $E_{n,\delta}(0)=0$ when $\delta\in[0,n+1)$ and $E_{n,n+1}(0)=1/(n+1)!$. We consider for $n=0,1,...$ and $\delta\in[n,n+1]$ the quantity
\beq \label{e20}
G=G_{n,\delta}(u)=({-}1)^{n+1}\,\frac{e^{-u}-\Bigl(1-\dfrac{u^1}{1!}+...+({-}1)^n \,\dfrac{u^n}{n!}\Bigr)}{u^{\delta}}~,~~~~u\geq0~,
\eq
with $G_{n,\delta}(0)=0$ when $\delta\in[n,n+1)$ and $G_{n,n+1}(0)=1/(n+1)!$. Then $E_{n,\delta}(s)$ and $G_{n,\delta}(u)$ are continuous functions of $s\geq0$ and $u\geq0$, respectively. We have for $n=0$ and $\delta\in[0,1]$
\beq \label{e21}
E_{0,\delta}(s)=\frac{1-e^{-s}}{s^{\delta}}\,,~~s\geq0~;~~~~~~ G_{0,\delta}(u)=\frac{1-e^{-u}}{u^{\delta}}\,,~~u\geq0~,
\eq
so that $E_{0,\delta}=G_{0,\delta}$. We have, see beginning of Secs.~\ref{sec3} and \ref{sec8}, for $n=0,1,...$
\beq \label{e22}
E_{n,\delta}(s)\geq0\,,~~s\geq0\,,~~\delta\in[0,n+1]~;~~~~~~ G_{n,\delta}(u)\geq0\,,~~u\geq0\,,~~\delta\in[n,n+1]~.
\eq

We have, see beginning of Sec.~\ref{sec3},
\beq \label{e23}
E_{n,0}(s)\leq1=\lim_{s\pr\infty}\,E_{n,0}(s)~;~~~~~~ E_{n,n+1}(s)\leq \frac{1}{(n+1)!}=\lim_{s\downarrow0}\,E_{n,n+1}(s)~.
\eq
We have, see beginning of Sec.~\ref{sec8},
\beq \label{e24}
G_{n,n}(u)\leq\frac{1}{n!}=\lim_{u\pr\infty}\,G_{n,n}(u)~;~~~~~~ G_{n,n+1}(u)\leq\frac{1}{(n+1)!}=\lim_{u\downarrow0}\,G_{n,n+1}(u)~.
\eq
Hence, the maximization problems for $E_{n,\delta}$ and $G_{n,\delta}$ have simple and explicit solutions for $\delta=0,n+1$ and $\delta=n,n+1$, respectively, and so we may confine attention to the cases $\delta\in(0,n+1)$ and $\delta\in(n,n+1)$, respectively.

We have chosen to develop the further results for the $S$-case and the $U$-case separately. While the type of results for both cases are globally the same, the precise formulation of the results, their proofs and the proof techniques used are often quite different. This separate treatment is meant to prevent reader's confusion that would occur when one has to switch back and forth between the $S$-case and $U$-case all the time. Accordingly, we describe the results for the $S$-case in Subsec.~\ref{subsec2.1} and those for the $U$-case in Subsec.~2.2. In Subsec.~\ref{subsec2.3} we summarize the further contents of this article and we point out where the various proofs of the results in Subsecs.~\ref{subsec2.1} and 2.2 can be found.

\subsection{Results for the $S$-case} \label{subsec2.1}
\mbox{} \\[-9mm]

We start with existence, uniqueness and characterization of stationary points $s>0$ of $E_{n,\delta}(s)$ when $n=0,1,...$ and $\delta\in(0,n+1)$.

\begin{prop} \label{prop1}
Let $n=0,1,..$ and $\delta\in(0,n+1)$. Then for $s>0$
\beq \label{e25}
\frac{d}{ds}\,[E_{n,\delta}(s)]=0\Leftrightarrow e^s=1+\frac{s^1}{1!}+...+ \frac{s^n}{n!}+\frac{1}{\delta}~\frac{s^{n+1}}{n!}~,
\eq
and the equation
\beq \label{e26}
e^s=1+\frac{s^1}{1!}+...+\frac{s^n}{n!}+\frac{1}{\delta}~\frac{s^{n+1}}{n!}
\eq
has exactly one positive solution $s>0$.
\end{prop}

\noindent
The unique $s>0$ such that (\ref{e26}) holds is denoted by $s_n(\delta)$. In the proof of Prop.~\ref{prop1} as given in Sec.~\ref{sec3}, we shall also establish that $s_n(\delta)>n+1-\delta$.

We thus have from Prop.~\ref{prop1} that
\beq \label{e27}
ME_{n,\delta}:=\max_{s\geq0}\,E_{n,\delta}(s)=E_{n,\delta}(s_n(\delta))~,
\eq
where we have also used that $E_{n,\delta}(s)\pr0$ as $s\downarrow0$ or $s\pr\infty$ when $n=0,1,...$ and $\delta\in(0,n+1)$.

\begin{prop} \label{prop2}
Let $n=0,1,...$ and $\delta\in(0,n+1)$. Then
\beq \label{e28}
ME_{n,\delta}=\frac{1-\Bigl(1+\dfrac{s^1}{1!}+...+\dfrac{s^n}{n!}\Bigr)\,e^{-s}} {s^{\delta}}=\frac{s^{n+1-\delta}}{n!\,\delta\,e^s}~,~~~~~~s=s_n(\delta)~.
\eq
\end{prop}

\begin{prop} \label{prop3}
Let $n=0,1,...$ and $\delta\in(0,n+1)$. Then
\beq \label{e29}
s_n'(\delta)=\frac{d}{d\delta}\,[s_n(\delta)]={-}\,\frac{1}{\delta}~\frac{s_n(\delta)} {s_n(\delta)-(n+1-\delta)}~.
\eq
\end{prop}

\noindent
We observe from $s_n(\delta)>n+1-\delta$ that $s_n'(\delta)<0$, and so $s_n(\delta)$ is strictly decreasing in $\delta\in(0,n+1)$.

It is of interest to find good bounds on $s_n(\delta)$, so that both the maximum location $s_n(\delta)$ and the maximum value $ME_{n,\delta}$ of $E_{n,\delta}(s)$, $s\geq0$, are accurately estimated. Given the context in which the $E_{n,\delta}$ arise, see the beginning of Sec.~\ref{sec1}, the need for sharpness of such bounds is in particular felt for the cases where $\delta$ is not much smaller than $n+1$. The bounds we derive do hold for all $\delta\in(0,n+1)$, but are indeed sharp in this preferred range with $\delta$ not far away from $n+1$.

A simple criterion to decide whether a particular $s>0$ is less than or exceeds $s_n(\delta)$ is as follows.

\begin{prop} \label{prop4}
Let $n=0,1,...$ and $\delta\in(0,n+1)$. Then for $s>0$
\beq \label{e30}
s<s_n(\delta)\Leftrightarrow e^s<1+\frac{s^1}{1!}+...+\frac{s^n}{n!}+\frac{1}{\delta}~ \frac{s^{n+1}}{n!}~,
\eq
\beq \label{e31}
s>s_n(\delta)\Leftrightarrow e^s>1+\frac{s^1}{1!}+...+\frac{s^n}{n!}+\frac{1}{\delta}~ \frac{s^{n+1}}{n!}~.
\eq
\end{prop}

\noindent
From Taylor's theorem with remainder in integral form, see \cite{ref3}, 1.4 (vi), applied to the function $f(s)=\exp(s)$, $s>0$, we have for $n=0,1,...$
\beq \label{e32}
e^s=1+\frac{s^1}{1!} +...+\frac{s^n}{n!}+\frac{1}{n!}\,e^s\,\il_0^s\,x^n\,e^{-x} \,dx~,~~~~~~s>0~.
\eq
Hence, we have from Prop.~\ref{prop4} for $n=0,1,...\,$, $\delta\in(0,n+1)$ and $s>0$
\beq \label{e33}
s<s_n(\delta)\Leftrightarrow\il_0^s\,x^n\,e^{-x}\,dx<\frac{1}{\delta}\,s^{n+1}\,e^{-s}~,
\eq
\beq \label{e34}
s>s_n(\delta)\Leftrightarrow\il_0^s\,x^n\,e^{-x}\,dx>\frac{1}{\delta}\,s^{n+1}\,e^{-s}~.
\eq
This observation is crucial in the proof of the following result.

\begin{prop} \label{prop5}
Let $n=0,1,...\,$.

a.~~Assume that $S(\delta)$ is a smooth function of $\delta\in(0,n+1]$ with $S(n+1)=0<S(\delta)$, $\delta\in(0,n+1)$, and that
\beq \label{e35}
S(\delta)>{-}\delta\,S'(\delta)(S(\delta)-(n+1-\delta))~,~~~~~~\delta\in(0,n+1)~.
\eq
Then $S(\delta)<s_n(\delta)$, $\delta\in(0,n+1)$.

b.~~Assume that $S(\delta)$ is a smooth function of $\delta\in(0,n+1]$ with $S(n+1)\geq0$, $S(\delta)>0$, $\delta\in(0,n+1)$, and that
\beq \label{e36} S(\delta)<{-}\delta\,S'(\delta)(S(\delta)-(n+1-\delta))~,~~~~~~\delta\in(0,n+1)~.
\eq
Then $S(\delta)>s_n(\delta)$, $\delta\in(0,n+1)$.
\end{prop}

\noindent
Observe that the $S(\delta)$ in Prop.~\ref{prop5}a and b satisfy (\ref{e29}), the differential equation for $s_n(\delta)$, with inequality instead of equality.

The inequality $s_n(\delta)>n+1-\delta$ is established at once from Prop.~\ref{prop5}a. There is the following sharpening.

\begin{prop} \label{prop6}
Let $n=0,1,...\,$. Then
\beq \label{e37}
{\rm ln}\Bigl(\frac{n+1}{\delta}\Bigr)+(n+1-\delta)<s_n(\delta)<(n+1-\delta)+(n+1-\delta)/\delta~, ~~~~~~\delta\in(0,n+1)~.
\eq
\end{prop}

\noindent
The two inequalities in (\ref{e37}) are reasonably sharp on the range $\delta\in[1,n+1)$, especially when $\delta$ is close to $n+1$. With $y=n+1-\delta$, we have the expansions
\beq \label{e38}
s_n(\delta)=\frac{n+2}{n+1}\,y+\frac{n+2}{n+3}~\frac{1}{(n+1)^2}\,y^2+...~,
\eq
\begin{eqnarray} \label{e39}
& \mbox{} & {\rm a.}~~{\rm ln}\Bigl(\frac{n+1}{\delta}\Bigr)+(n+1-\delta)=\frac{n+2}{n+1}\,y+\frac12~ \frac{1}{(n+1)^2}\,y^2+...~, \nonumber \\[3mm]
& & {\rm b.}~~(n+1-\delta)+(n+1-\delta)/\delta= \frac{n+2}{n+1}\,y+\frac{1}{(n+1)^2}\,y^2+...~.
\end{eqnarray}
From the two inequalities in Prop.~\ref{prop6} together with Prop.~\ref{prop3}, it is seen that $s_n(\delta)$ decreases strictly from $+\infty$ at $\delta=0$ to 0 at $\delta=n+1$. Furthermore, there is the following result.

\begin{prop} \label{prop7}
Let $n=0,1,...\,$. Then $s_n(\delta)$ is strictly convex in $\delta\in(0,n+1)$.
\end{prop}

\noindent
The proof of Prop.~\ref{prop7} consists of manipulating the convexity condition $s_n''(\delta)>0$ by using Prop.~\ref{prop3} to the equivalent condition
\beq \label{e38a}
s_n(\delta)>y+(y+\tfrac14\,\delta^2)^{1/2}-\tfrac12\,\delta~,~~~~~~ y=n+1-\delta\in(0,n+1)~.
\eq
The inequality in (\ref{e38a}) can be established from Prop.~\ref{prop5}a, see the proof of Lemma~\ref{lem1} in Sec.~\ref{sec4}.

While the lower bound in Prop.~\ref{prop6} is reasonably sharp, the upper bound is not. For instance, the asymptotics of $s_0(\delta)$ as $\delta\downarrow0$ is given by, see Lemma~\ref{lem3},
\beq \label{e39a}
s_0(\delta)=z+{\rm ln}\,z+O\Bigl(\frac{{\rm ln}\,z}{z}\Bigr)~,~~~~~~ z={\rm ln}\Bigl(\frac{1}{\delta}\Bigr)\,,~~\delta\downarrow0~.
\eq
A $(z+{\rm ln}\,z)$-behaviour for general $n=0,1,...$ is captured by the following result.

\begin{prop} \label{prop8}
Let $n=0,1,...\,$. Then we have for $\delta\in(0,n+1)$
\beq \label{e40}
s_n(\delta)<z+(n+1-\delta)\,{\rm ln}\,z-(n+1)\,\frac{\partial\gamma}{\partial a}\,(1,z)~,~~~~~~ z={\rm ln}\Bigl(\frac{n+1}{\delta}\Bigr)~,
\eq
where $\gamma(a,z)$ is the incomplete Gamma function,
\beq \label{e41}
\gamma(a,z)=\il_0^z\,t^{a-1}\,e^{-t}\,dt~,~~~~~~z>0\,,~~a>0~,
\eq
see \cite{ref3}, Ch.~8, so that
\beq \label{e42}
\frac{\partial\gamma}{\partial a}\,(1,z)=\il_0^z\,e^{-t}\,{\rm ln}\,t\,dt~,~~~~~~z>0~.
\eq
\end{prop}

\noindent
Observe that the right-hand side of (\ref{e39a}) is singular at $\delta=1$, $z=0$, while the right-hand side of (\ref{e40}) tends to 0 as $\delta\uparrow n+1$. The right-hand side of (\ref{e40}) is somewhat awkward to deal with when trying to use Prop.~\ref{prop5}b. Therefore, the proof of Prop.~\ref{prop8} uses a different approach in which we insert the first inequality in (\ref{e37}) into the expression for $s_n'(\delta)$ as given by Prop.~\ref{prop3} to obtain a lower bound for $s_n'(\delta)$, $\delta\in(0,n+1)$. Using that $s_n(n+1)=0$, this lower bound for $s_n'(\delta)$ can be converted to an upper bound for $s_n(\delta)$ by integrating it from $\delta$ to $n+1$.

We next consider the function, see Prop.~\ref{prop2},
\beq \label{e43}
ME_{n,\delta}=E_{n,\delta}(s_n(\delta))=\frac{s^{n+1-\delta}}{\delta\cdot n!\,e^s}~, ~~~~~~s=s_n(\delta)~.
\eq

\begin{prop} \label{prop9}
Let $n=0,1,...\,$. Then
\beq \label{e44}
ME_{n,\delta}\pr1\,,~~\delta\downarrow0~;~~~~~~ME_{n,\delta}\pr \frac{1}{(n+1)!}\,,~~\delta\uparrow n+1~.
\eq
\end{prop}

\noindent
Thus, see (\ref{e23}), we have that $ME_{n,\delta}$ is continuous at $\delta=0$ and $\delta=n+1$. The proof of Prop.~\ref{prop9} as given in Sec.~\ref{sec4} uses Props.~\ref{prop6} and \ref{prop8}.

\begin{prop} \label{prop10}
Let $n=0,1,...\,$. Then
\beq \label{e45}
\frac{d}{d\delta}\,[{\rm ln}\,ME_{n,\delta}]={-}{\rm ln}(s_n(\delta))~, ~~~~~~\delta\in(0,n+1)~.
\eq
\end{prop}

\noindent
The proof of Prop.~\ref{prop10} consists of a computation, using the concise form of $ME_{n,\delta}$ in (\ref{e43}) and the expression for $s_n'(\delta)$ as given in Prop.~\ref{prop3}. From Prop.~\ref{prop10} and $s_n'(\delta)<0<s_n(\delta)$, $\delta\in(0,n+1)$, it follows that $ME_{n,\delta}$ is a log-convex function of $\delta\in(0,n+1)$. In particular, $ME_{n,\delta}$ is a convex function of $\delta\in(0,n+1)$.

While $s_n(\delta)$, and therefore $ME_{n,\delta}=E_{n,\delta}(s_n(\delta))$, must be computed numerically, see Sec.~\ref{sec3} for this issue, we can give a simple, closed-form expression for the minimum of $ME_{n,\delta}$ over $\delta\in(0,n+1)$.

\begin{prop} \label{prop11}
Let $n=0,1,...\,$. Then the minimum $\hat{E}_n$ of $ME_{n,\delta}$ over $\delta\in(0,n+1)$ is given by
\beq \label{e46}
\hat{E}_n=\frac{1}{e}\,\Bigl(e-\sum_{k=0}^n\,\frac{1}{k!}\Bigr)~,
\eq
and is assumed for
\beq \label{e47}
\delta=\delta_{n,E}=\Bigl(n!\Bigl(e-\sum_{k=0}^n\,\frac{1}{k!}\Bigr)\Bigr)^{-1}~.
\eq
\end{prop}

\noindent
The proof of Prop.~\ref{prop11} consists of the observation that $ME_{n,\delta}$ is minimal for the unique $\delta\in(0,n+1)$ such that $s_n(\delta)=1$, see Prop.~\ref{prop10}. This is combined with (\ref{e26}) and (\ref{e28}), used with this $\delta$.

When we have an approximation $S(\delta)$ of $s_n(\delta)$, then $E_{n,\delta}(S(\delta))$ is a lower bound for $ME_{n,\delta}$. We describe now a method to find an upper bound for $ME_{n,\delta}$ given a lower bound $S(\delta)$ for $s_n(\delta)$ satisfying $S(\delta)\geq n+1-\delta$, $\delta\in(0,n+1)$. To that end, we define for $s>0$
\beq \label{e48}
F_1(s)=\frac{s^{n+1-\delta}}{\delta\cdot n!\,e^s}~,~~~~~~ F_2(s)=\frac{1}{\delta\cdot n!}~\frac{s^{n+1-\delta}} {1+\dfrac{s^1}{1!}+...+\dfrac{s^n}{n!}+\dfrac{1}{\delta}~\dfrac{s^{n+1}}{n!}}~,
\eq
so that, see (\ref{e28}),
\beq \label{e49}
ME_{n,\delta}=F_1(s_n(\delta))=F_2(s_n(\delta))~,~~~~~~ \delta\in(0,n+1)~.
\eq

\begin{prop} \label{prop12}
Both $F_1$ and $F_2$ are strictly decreasing in $s\in(n+1-\delta,s_n(\delta))$. Moreover, $F_2(s)<F_1(s)$, $s\in(n+1-\delta,s_n(\delta))$.
\end{prop}

\noindent
When now $S(\delta)\in[n+1-\delta,s_n(\delta))$, we have from (\ref{e49}) and Prop.~\ref{prop12} for $i=1,2$
\beq \label{e50}
E_{n,\delta}(s_n(\delta))=F_i(s_n(\delta))<F_i(S(\delta))~,~~~~~~ \delta\in(0,n+1)~,
\eq
while $F_2(S(\delta))<F_1(S(\delta))$. Hence, both $F_1(S(\delta))$ and $F_2(S(\delta))$ are an upper bound for $ME_{n,\delta}$ and $F_2(S(\delta))$ is the sharpest of the two.

\subsection{Results for the $U$-case} \label{subsec2.2}
\mbox{} \\[-9mm]

We consider $n=1,2,...\,$, the case $n=0$ already being covered by the results of Subsec.~\ref{subsec2.1}. We first note that $G_{n,\delta}(u)>0$ for $u>0$, where $G_{n,\delta}$ is as in (\ref{e20}). Indeed, when $u>0$ there is by Taylor's theorem a $\xi\in(0,u)$ such that
\beq \label{e51}
e^{-u}-\Bigl(1-\frac{u^1}{1!}+...+({-}1)^n\,\frac{u^n}{n!}\Bigr)=({-}1)^{n+1}\, \frac{u^{n+1}}{(n+1)!}\,e^{-\xi}~.
\eq
We next consider existence, uniqueness and characterization of stationary points $u>0$ of $G_{n,\delta}(u)$.

\begin{prop} \label{prop13}
Let $n=1,2,...$ and $\delta\in(n,n+1)$. Then for $u>0$
\beq \label{e52}
\frac{d}{du}\,[G_{n,\delta}(u)]=0\Leftrightarrow e^{-u}=1-\frac{u^1}{1!} +...+({-}1)^{n-1}\,\frac{u^{n-1}}{(n{-}1)!}+({-}1)^n\,\frac{\delta u^n} {n!\,(u{+}\delta)}\,,
\eq
and the equation
\beq \label{e53}
e^{-u}=1-\frac{u^1}{1!}+...+({-}1)^{n-1}\,\frac{u^{n-1}}{(n-1)!}+({-}1)^n \,\frac{\delta\,u^n}{n!\,(u+\delta)}
\eq
has exactly one positive solution $u$.
\end{prop}

\noindent
The unique $u>0$ such that (\ref{e53}) holds is denoted by $u_n(\delta)$. As a consequence of the proof of Prop.~\ref{prop13} we have that $u_n(\delta)$ is a strictly decreasing function of $\delta\in (n,n+1)$.

We thus have from Prop.~\ref{prop13} that
\begin{eqnarray} \label{e54}
& \mbox{} & MG_{n,\delta}:= \max_{u\geq0}\,G_{n,\delta}(u)=G_{n,\delta}(u_n(\delta)) \nonumber \\[3mm]
& & =~ ({-}1)^{n+1}\,\frac{e^{-u}-\Bigl(1-\dfrac{u^1}{1!}+...+({-}1)^n\,\dfrac{u^n}{n!}\Bigr)} {u^{\delta}}=\frac{1}{n!}~\frac{u^{n+1-\delta}}{u+\delta}~.
\end{eqnarray}
Here we have also used that $G_{n,\delta}(u)\pr0$ as $u\downarrow0$ or $u\pr\infty$ when $n=1,2,...$ and $\delta\in(n,n+1)$. Furthermore, in the second line of (\ref{e54}) we have taken $u=u_n(\delta)$; the last identity in (\ref{e54}) is proved at the end of Sec.~\ref{sec8}, see Lemma~\ref{lem5}.

\begin{prop} \label{prop14}
Let $n=1,2,...$ and $\delta\in(n,n+1)$. Then
\beq \label{e55}
u_n'(\delta)=\frac{d}{d\delta}\,[u_n(\delta)]={-}\,\frac{u_n(\delta)} {(\delta-n)\,u_n(\delta)-(n+1-\delta)\,\delta}~.
\eq
\end{prop}

\noindent
It is observed in the proof of Prop.~\ref{prop14} as given in Sec.~\ref{sec8} that $(\delta-n)\,u_n(\delta)-(n+1-\delta)\,\delta>0$ for $\delta\in(n,n+1)$. A stronger inequality follows from Prop.~\ref{prop17}, the proof of which does not use Prop.~\ref{prop14}.

We now turn to bounds on $u_n(\delta)$. A simple criterion to decide whether a particular $u>0$ is less than or exceeds $u_n(\delta)$ is as follows.

\begin{prop} \label{prop15}
Let $n=1,2,...$ and $\delta\in(n,n+1)$. Then for $u>0$
\begin{eqnarray} \label{e56}
& \mbox{} & \hspace*{-7mm}u<u_n(\delta) \nonumber \\[3mm]
& & \hspace*{-7mm}\Leftrightarrow~ ({-}1)^n\Bigl(e^{-u}-\Bigl(1-\frac{u^1}{1!}+...+({-}1)^{n-1}\, \frac{u^{n-1}}{(n-1)!}+({-}1)^n\,\frac{\delta u^n}{n!(u+\delta)}\Bigr)\Bigr)>0~, \nonumber \\
\mbox{}
\end{eqnarray}
\begin{eqnarray} \label{e57}
& \mbox{} & \hspace*{-7mm}u>u_n(\delta) \nonumber \\[3mm]
& & \hspace*{-7mm}\Leftrightarrow~ ({-}1)^n\Bigl(e^{-u}-\Bigl(1-\frac{u^1}{1!}+...+({-}1)^{n-1}\, \frac{u^{n-1}}{(n-1)!}+({-}1)^n\,\frac{\delta u^n}{n!(u+\delta)}\Bigr)\Bigr)<0~, \nonumber \\
\mbox{}
\end{eqnarray}
\end{prop}
\mbox{} \\
From Taylor's formula with remainder in integral form, see \cite{ref3}, 1.4 (vi), applied to the function $f(u)=\exp({-}u)$, $u>0$, we have for $n=1,2,...$
\beq \label{e58}
e^{-u}=1-\frac{u^1}{1!}+...+({-}1)^{n-1}\,\frac{u^{n-1}}{(n-1)!}+ \frac{({-}1)^n}{(n-1)!}\,e^{-u}\,\il_0^u\,x^{n-1}\,e^x\,dx~.
\eq
Hence, we have for $n=1,2,...\,$, $\delta\in(n,n+1)$ and $u>0$ that
\beq \label{e59}
u<u_n(\delta)\Leftrightarrow\il_0^u\,x^{n-1}\,e^x\,dx>\frac1n~\frac{\delta u^n} {u+\delta}\,e^u~,
\eq
\beq \label{e60}
u>u_n(\delta)\Leftrightarrow\il_0^u\,x^{n-1}\,e^x\,dx<\frac1n~\frac{\delta u^n} {u+\delta}\,e^u~.
\eq
This observation leads to the following result.

\begin{prop} \label{prop16}
Let $n=1,2,...\,$.

a.~~Assume that $U(\delta)$ is a smooth function of $\delta\in(n,n+1]$ with $U(n+1)=0<U(\delta)$, $\delta\in(n,n+1)$, and that
\beq \label{e61}
U(\delta)>{-}U'(\delta)((\delta-n)\,U(\delta)-(n+1-\delta)\,\delta)~, ~~~~~~\delta\in(n,n+1)~.
\eq
Then $U(\delta)<u_n(\delta)$, $\delta\in(n,n+1)$.

b.~~Assume that $U(\delta)$ is a smooth function of $\delta\in(n,n+1]$ with $U(n+1)\geq0$, $U(\delta)>0$, $\delta\in(n,n+1)$, and that
\beq \label{e62}
U(\delta)<{-}U'(\delta)((\delta-n)\,U(\delta)-(n+1-\delta)\,\delta)~,~~~~~~ \delta\in(n,n+1)~.
\eq
Then $U(\delta)>u_n(\delta)$, $\delta\in(n,n+1)$.
\end{prop}

\noindent
Observe that the $U(\delta)$ in Prop.~\ref{prop16}a and b satisfy (\ref{e55}), the differential equation for $u_n(\delta)$, with inequality instead of equality.

The inequality $u_n(\delta)>(n+1-\delta)\,\delta/(\delta-n)$, $\delta\in(n,n+1)$ is established at once from Prop.~\ref{prop16}a. There is the following sharpening.

\begin{prop} \label{prop17}
Let $n=1,2,...\,$. Then
\beq \label{e63}
\frac{(n+1-\delta)\,\delta}{\delta-n}+(n+1-\delta)<u_n(\delta)<\frac{\delta} {\delta-n}~,~~~~~~\delta\in(n,n+1)~.
\eq
\end{prop}

\noindent
The two inequalities in (\ref{e63}) are reasonably sharp when $\delta\downarrow n$:
\beq \label{e64}
\frac{\delta}{\delta-n}-\Bigl(\frac{(n+1-\delta)\,\delta}{\delta-n}+ (n+1-\delta)\Bigr)=\delta-(n+1-\delta)<\delta~,~~~~~~\delta\in(n,n+1)~,
\eq
so that the difference between upper bound and lower bound in (\ref{e63}) remains bounded as $\delta\downarrow n$. The upper bound in (\ref{e63}) is not sharp as $\delta\uparrow n+1$.

The following results show that there are upper bounds for $u_n(\delta)$ that are sharp, both as $\delta\downarrow n$ and $\delta\uparrow(n+1)$.

\begin{prop} \label{prop18}
Let $n=2,3,...\,$. Then
\beq \label{e65}
u_n(\delta)<\frac{(n+1-\delta)\,\delta}{\delta-n}+\min\,\{1,n+1-\delta+\tfrac12\,(n+1-\delta)^2\}~, ~~~~~~\delta\in(n,n+1)~.
\eq
\end{prop}

\noindent
Observe that $n+1-\delta+\frac12\,(n+1-\delta)^2<1$ when $n+2-\sqrt{3}<\delta<n+1$.

\begin{prop} \label{prop19}
We have
\beq \label{e66}
u_1(\delta)<\frac{(2-\delta)\,\delta}{\delta-1}+(2-
\delta)+\tfrac12\,(2-\delta)^2+\tfrac13\,(2-\delta)^3~,~~~~~~\delta\in(1,2)~.
\eq
\end{prop}

\noindent
From Props.~\ref{prop18} and \ref{prop19} there follows for $n=1,2,...$ the simple and convenient bound
\beq \label{e67}
u_n(\delta)<\frac{(n+1-\delta)\,\delta}{\delta-n}-{\rm ln}(\delta-n)~,~~~~~~ \delta\in(n,n+1)~.
\eq
Observe also that the difference between the upper bound in (\ref{e65}) for $u_n(\delta)$ and the lower bound in (\ref{e63}) for $u_n(\delta)$ equals $\min\,\{\delta-n,\frac12\,(n+1-\delta)^2\}\pr0$ as $\delta\downarrow n$ or $\delta\uparrow(n+1)$, case $n=2,3,...\,$. Furthermore, for $n=2,3,...\,$, we have from the upper bound in (\ref{e65}) for $u_n(\delta)$
\beq \label{e68}
\frac{\delta}{\delta-n}-u_n(\delta)>\delta-1>n-1>0~,~~~~~~\delta\in(n,n+1)~.
\eq
Hence, for $n=2,3,...\,$, the upper bound in (\ref{e65}) for $u_n(\delta)$ is truly sharper than the upper bound in (\ref{e63}). The situation for $n=1$ is quite different.

\begin{prop} \label{prop20}
We have
\beq \label{e69}
u_1(\delta)=B(1+O(Be^{-B}))~,~~~~~~B=\frac{\delta}{\delta-1}\,,~~ \delta\downarrow1~.
\eq
\end{prop}

\noindent
Hence, $u_1(\delta)-\delta/(\delta-1)$ is exponentially small as $\delta\downarrow1$.

It is seen from Props.~\ref{prop14}, \ref{prop17}, \ref{prop18} and \ref{prop19} that $u_n(\delta)$ decreases strictly from $+\infty$ at $\delta=n$ to 0 at $\delta=n+1$, $n=1,2,...\,$. Furthermore, there is the following result.

\begin{prop} \label{prop21}
Let $n=1,2,...\,$. Then $u_n(\delta)$ is strictly convex in $\delta\in(n,n+1)$.
\end{prop}

\noindent
The proof of Prop.~\ref{prop21}, as given in Sec.~\ref{sec9}, consists of combining the expression in Prop.~\ref{prop14} for $u_n'(\delta)$ and the first inequality in (\ref{e63}), which leads to the inequality
\beq \label{e70}
u_n'(\delta)>{-}\,\frac{2\delta-n}{(\delta-n)^2}~,~~~~~~\delta\in(n,n+1)~.
\eq
The convexity condition $u_n''(\delta)>0$, $\delta\in(n,n+1)$ can be written as
\beq \label{e71}
u_n^2(\delta)+(n+1-\delta)\,\delta\,u_n'(\delta)+(2\delta-n-1)\,u_n(\delta)>0~, ~~~~~~\delta\in(n,n+1)~.
\eq
Using (\ref{e70}) and the first inequality in (\ref{e63}), the inequality in (\ref{e71}) readily follows.

We next consider the function, see (\ref{e54}),
\beq \label{e72}
MG_{n,\delta}=G_{n,\delta}(u_n(\delta))=\frac{1}{n!}~\frac{(u_n(\delta))^
{n+1-\delta}}{u_n(\delta)+\delta}~,~~~~~~\delta\in(n,n+1)~.
\eq

\begin{prop} \label{prop22}
Let $n=1,2,...\,$. Then
\beq \label{e73}
MG_{n,\delta}\pr\frac{1}{n!}\,,~~\delta\downarrow n~;~~~~~~ MG_{n,\delta}\pr \frac{1}{(n+1)!}\,,~~\delta\uparrow(n+1)~.
\eq
\end{prop}

\noindent
Thus, see (\ref{e24}), we have that $MG_{n,\delta}$ is continuous at $\delta=n$ and $\delta=n+1$. The proof of Prop.~\ref{prop22} uses (\ref{e67}) and Prop.~\ref{prop17}.

\begin{prop} \label{prop23}
Let $n=1,2,...\,$. Then
\beq \label{e74}
\frac{d}{d\delta}\,[{\rm ln}\,MG_{n,\delta}]={-}{\rm ln}(u_n(\delta))~, ~~~~~~\delta\in(n,n+1)~.
\eq
\end{prop}

\noindent
The proof of Prop.~\ref{prop23} consists of a computation in which the expression for $u_n'(\delta)$ as given in Prop.~\ref{prop14} and the concise form for $MG_{n,\delta}$ in (\ref{e72}) are used.

Since $u_n'(\delta)<0<u_n(\delta)$, it is seen that $MG_{n,\delta}$ is a log-convex function of $\delta\in(n,n+1)$. Hence, $MG_{n,\delta}$ is a convex function of $\delta\in(n,n+1)$. The following result shows that the minimum of $MG_{n,\delta}$ over $\delta\in(n,n+1)$ has a simple form.

\begin{prop} \label{prop24}
Let $n=1,2,...\,$. Then the minimum $\hat{G}_n$ of $MG_{n,\delta}$ over $\delta\in(n,n+1)$ is given by
\beq \label{e75}
\hat{G}_n=({-}1)^{n+1}\,\Bigl(e^{-1}-\Bigl(1-\frac{1}{1!}+...+ ({-}1)^n\,\frac{1}{n!}\Bigr)\Bigr)~,
\eq
and is assumed for
\beq \label{e76}
\delta=\delta_{n,G}=\Bigl[({-}1)^{n+1}\,n!\Bigl( e^{-1}-\Bigl(1-\frac{1}{1!} +...+({-}1)^n\,\frac{1}{n!}\Bigr)\Bigr)\Bigr]^{-1}-1~.
\eq
\end{prop}

\noindent
The proof of Prop.~\ref{prop24} consists of the observation that $MG_{n,\delta}$ is minimal for the unique $\delta\in(n,n+1)$ such that $u_n(\delta)=1$, see Prop.~\ref{prop23}. This is combined with (\ref{e20}) and (\ref{e72}) with this $\delta$.

The approximations $U(\delta)$ of $u_n(\delta)$ given in Props.~\ref{prop17}, \ref{prop18} and \ref{prop19} give rise to various lower bounds for $MG_{n,\delta}$ in the form $G_{n,\delta}(U(\delta))$. Now let $U(\delta)$ be a lower bound for $u_n(\delta)$ that satisfies $U(\delta)\geq(n+1-\delta)\,\delta/(\delta-n)$, $\delta\in(n,n+1)$, and define
\beq \label{e77}
H(u)=\frac{1}{n!}~\frac{u^{n+1-\delta}}{u+\delta}~,~~~~~~u\geq0~.
\eq
Then $H(u)$ is strictly decreasing in $u\in((n+1-\delta)\,\delta/(\delta-n),u_n(\delta))$. Therefore, by (\ref{e72}),
\beq \label{e78}
MG_{n,\delta}=G_{n,\delta}(u_n(\delta))=H(u_n(\delta))<H(U(\delta))~.
\eq
Hence $H(U(\delta))$ is an upper bound for $MG_{n,\delta}$.

\subsection{Organization of the remainder of this article} \label{subsec2.3}
\mbox{} \\[-9mm]

An examination of the results presented in Subsecs.~\ref{subsec2.1} and 2.2 for the $S$-case and $U$-case, respectively, shows that there is a common plan according to which these results emerge. There are, however, obvious differences in the formulation, setting, tightness and proof details. We list these below.
\bi{--0}
\ITEM{--} the functions $s_n(\delta)$ have $\delta\in(0,n+1]$ while the functions $u_n(\delta)$ have $\delta\in(n,n+1]$,
\ITEM{--} the results for $s_n$ hold for all $n=0,1,...$ while the results for $u_n$ hold for $n=1,2,...\,$, with the case $n=1$ still being somewhat exceptional,
\ITEM{--} the bounds for $u_n$ are clearly sharper than those for $s_n$,
\ITEM{--} the $u_n$ exhibit a $((n+1-\delta)\,\delta/(\delta-n)+(n+1-\delta))$-behaviour on the whole range $\delta\in(n,n+1)$ while the $s_n$ exhibit a similar behaviour on much of the range $\delta\in(0,n+1)$ but with a particular $(z+{\rm ln}\,z)$-behaviour, $z={\rm ln}((n+1)/\delta)$ as $\delta\downarrow0$,
\ITEM{--} the approach based on Taylor's formula with remainder in integral form to find bounds, finds more fertile ground in the $U$-case than that it does in the $S$-case,
\ITEM{--} an additional technique, based on the differential equation satisfied by the $s_n$, is required to bring out the particular $(z+{\rm ln}\,z)$-behaviour as $\delta\downarrow0$; such an effort is not required for the $u_n$,
\ITEM{--} the minimal values $\hat{E}_n$ and $\hat{G}_n$ of $ME_{n,\delta}$, $\delta\in(0,n+1)$ and $MG_{n,\delta}$, $\delta\in(n,n+1)$, respectively, are assumed for a $\delta=\delta_{n,E}$ and a $\delta=\delta_{n,G}$ that are relatively close to the endpoint $\delta=n+1$; while $\hat{G}_n$ is just marginally smaller than the value $1/(n+1)!$ of $MG_{n,\delta}$ at $\delta=n+1$, the $\hat{E}_n$ is about a factor $e$ smaller than the value $1/(n+1)!$ of $ME_{n,\delta}$ at $\delta=n+1$.
\ei

The circumstance that the results for the $S$-case and the $U$-case exhibit, on zooming in, such clear differences, leads us to develop and prove these results separately. As a consequence, the results have been presented into the separate Subsec.~\ref{subsec2.1} and 2.2, and this then also determines the organization of the remainder of this article. In Sec.~~\ref{sec3} we present basic facts about $E_{n,\delta}$ and $s_n$, such as existence uniqueness and (numerical) computation of the $s_n$. The computation of the $s_n$ can be done, in principle, by using Newton's method to solve the defining equation (\ref{e26}) for $s_n(\delta)$. However, when $n$ is somewhat larger and $\delta$ is close to $n+1$, this becomes awkward, and then a series representation of $s_n(\delta)$ in powers of $(n+1-\delta)$, based on the B\"urmann-Lagrange inversion formula is more appropriate. In Sec.~\ref{sec4} we present the bounds for $s_n$. This includes the bounds, and some of their consequences, that can be obtained using the approach based on Taylor's formula with remainder in integral form, as well as the bound in Prop.~\ref{prop8} that is obtained by exploring the expression in Prop.~\ref{prop3} for $s_n'(\delta)$. In Sec.~\ref{sec5} we show log-convexity of $ME_{n,\delta}$ as a function of $\delta\in(0,n+1)$, and we determine the minimum value of $ME_{n,\delta}$ as well as its minimizer $\delta=\delta_{n,E}$ in closed form. In Sec.~\ref{sec6} we detail the technique, explained at the end of Subsec.~\ref{subsec2.1} to pass from a lower bound $S(\delta)$ to an upper bound for $ME_{n,\delta}$, $\delta\in(0,n+1)$. In Sec.~\ref{sec7} we illustrate the results of Secs.~\ref{sec4}, \ref{sec5} and \ref{sec6} for $s_n$ and $ME_{n,\delta}$ with $n=1$ and $n=3$.

In Sec.~\ref{sec8} we present basic facts about $G_{n,\delta}$ and $u_n$, in a similar manner as we did in Sec.~\ref{sec3} for $E_{n,\delta}$ and $s_n$, although the proof details can be quite different at places. In Sec.~\ref{sec9} we present the bounds for $u_n$, and some of their consequences. These bounds are obtained by using the approach based on Taylor's formula with remainder in integral form, and are already so sharp that additional techniques, such as the one in Sec.~\ref{sec4} for $s_n$, are not required. We pay special attention in Sec.~\ref{sec9} to the somewhat exceptional case $n=1$. In Sec.~\ref{sec10} we show log-convexity of $MG_{n,\delta}$ as a function of $\delta\in(n,n+1)$, and we determine the minimum value of $MG_{n,\delta}$ as well as its minimizer $\delta=\delta_{n,G}$ in closed form. In Sec.~\ref{sec11} we illustrate the results of Secs.~\ref{sec9} and \ref{sec10} and of the technique given at the end of Subsec.~2.2 to pass from a lower bound $U(\delta)$ for $u_n(\delta)$ to an upper bound for $MG_{n,\delta}$, $\delta\in(n,n+1)$, for $u_n$ and $MG_{n,\delta}$ with $n=1$ and $n=3$.

\section{Basic properties of $E_{n,\delta}$ and $s_n$} \label{sec3}
\mbox{} \\[-9mm]

We recall that $E_{n,\delta}$ is defined for $n=0,1,...\,$, $\delta\in[0,n+1]$ and $s>0$ by
\beq \label{e79}
E_{n,\delta}(s)=\frac{1-\Bigl(1+\dfrac{s^1}{1!}+...+\dfrac{s^n}{n!}\Bigr)\,e^{-s}} {s^{\delta}}~,
\eq
and that
\beq \label{e80}
E_{n,\delta}(0)=\dlim_{s\downarrow0}\,E_{n,\delta}(s)= \left\{\ba{llll}
0 & \!\!, & ~~~\delta\in[0,n+1) & \!\!, \\[3mm]
\dfrac{1}{(n+1)!} & \!\!, & ~~~\delta=n+1 & \!\!.
\ea\right.
\eq
Then $E_{n,\delta}(s)$ is a continuous function of $s\geq0$. For the case $n=0$, $\delta\in[0,1]$, we have
\beq \label{81}
E_{0,\delta}(s)=\frac{1-e^{-s}}{s^{\delta}}\,,~~s>0~;~~~~~~ E_{0,\delta}(0)= \left\{\ba{llll}
0 & \!\!, & ~~~\delta\in[0,1) & \!\!, \\[2mm]
1 & \!\!, & ~~~\delta=1 & \!\!.
\ea\right.
\eq
We begin by establishing (\ref{e23}). \\ \\
{\bf Proof of (\ref{e23}).}~~We have for $n=0,1,...$ and $\delta=0$
\beq \label{e82}
0\leq E_{n,0}(s)=1-\Bigl(1+\frac{s^1}{1!}+...+\frac{s^n}{n!}\Bigr)\,e^{-s}\uparrow 1~,~~~~~~s\pr\infty~,
\eq
where monotonicity follows from $\frac{d}{ds}\,[E_{n,0}(s)]=\frac{1}{n!}\,s^n\,e^{-s}>0$.

We have for $n=0,1,...$ and $\delta=n+1$
\begin{eqnarray} \label{e83}
& \mbox{} & \hspace*{-8mm}0\leq E_{n,n+1}(s)=\frac{1}{(n+1)!}~\frac{1+\dfrac{1}{n+2}\,s^1+\dfrac{1} {(n+1)(n+2)}\,s^2+...}{1+\dfrac{1}{1!}\,s^1+\dfrac{1}{2!}\,s^2+...}\uparrow \frac{1}{(n+1)!}~, \nonumber \\
& & \hspace*{10.5cm}s\downarrow0~.
\end{eqnarray}
Thus, both items in (\ref{e23}) hold.\hfill \qedsymbol{} \\ \\
We infer from (\ref{e82}) and (\ref{e83}) that
\beq \label{e84}
\max_{s\geq0}\,E_{n,0}(s)=1~;~~~~~~\max_{s\geq0}\,E_{n,n+1}(s)=\frac{1}{(n+1)!}~.
\eq
{\bf Proof of Proposition~\ref{prop1}.}~~The equivalence in (\ref{e25}) follows from the definition of $E_{n,\delta}(s)$ in (\ref{e79}) with $n=0,1,...\,$, $\delta\in(0,n+1)$ and $s>0$ and a straightforward calculation.

Next, since $0<\delta<n+1$, we have
\begin{eqnarray} \label{e85}
e^s & = & 1+\frac{s^1}{1!}+...+\frac{s^n}{n!}+\frac{s^{n+1}}{(n+1)!}+... \nonumber \\[3mm]
& < & 1+\frac{s^1}{1!}+...+\frac{s^n}{n!}+\frac{1}{\delta}~\frac{s^{n+1}}{n!}
\end{eqnarray}
for small positive $s$. Obviously,
\beq \label{e86}
e^s>1+\frac{s^1}{1!}+...+\frac{s^n}{n!}+\frac{1}{\delta}~\frac{s^{n+1}}{n!}
\eq
for large positive $s>0$. The equation in (\ref{e26}) can be written as
\beq \label{e87}
\psi(s):=s-{\rm ln}\Bigl(1+\frac{s^1}{1!}+...+\frac{s^n}{n!}+\frac{1}{\delta}~ \frac{s^{n+1}}{n!}\Bigr)=0~.
\eq
We compute for $s>0$
\beq \label{e88}
\psi'(s)=\frac{1}{\delta\cdot n!}~\frac{s^n(s-(n+1-\delta))} {1+\dfrac{s^1}{1!}+...+\dfrac{s^n}{n!}+\dfrac{s^{n+1}}{\delta\cdot n!}}~.
\eq
Hence, $\psi$ is negative and strictly decreasing in $s\in(0,n+1-\delta)$, and $\psi$ is strictly increasing in $s>n+1-\delta$. Since $\psi(s)$ is positive for large positive $s$, we conclude that $\psi$ has exactly one positive zero which exceeds $n+1-\delta$.\hfill \qedsymbol{} \\ \\
{\bf Proof of Proposition~\ref{prop4}.}~~With $s_n(\delta)$ the unique, positive zero of $\psi$ where $n=0,1,...$ and $\delta\in(0,n+1)$, we have from the proof of Prop.~\ref{prop1} that $0<s<s_n(\delta)\Leftrightarrow\psi(s)<0$ and $s>s_n(\delta)\Leftrightarrow\psi(s)>0$.\hfill \qedsymbol{} \\ \\
{\bf Proof of Proposition~\ref{prop2}.}~~This follows at once from the defining equation (\ref{e26}) for $s_n(\delta)$.\hfill \qedsymbol{} \\ \\
{\bf Proof of Proposition~\ref{prop3}.}~~Writing $s=s_n(\delta)$ and $s'=s_n'(\delta)=\frac{d}{ds}\,[s_n(\delta)]$, we get by implicit differentiation in (\ref{e26})
\beq \label{e89}
e^s\,s'=\Bigl(1+\frac{s^1}{1!}+...+\frac{s^{n-1}}{(n-1)!}\Bigr)\,s'-\frac{1}{\delta^2}~ \frac{s^{n+1}}{n!}+\frac{1}{\delta}~\frac{n+1}{n!}\,s^n\,s'~,
\eq
i.e., that
\beq \label{e90}
\Bigl(e^s-1-\frac{s^1}{1!}-...-\frac{s^{n-1}}{(n-1)!}-\frac{1}{\delta}~ \frac{n+1}{n!}\,s^n\Bigr)\,s'={-}\,\frac{1}{\delta^2}~\frac{s^{n+1}}{n!}~.
\eq
Using (\ref{e26}) this gives
\beq \label{e91}
\Bigl(\frac{s^n}{n!}+\frac{1}{\delta}~\frac{s^{n+1}}{n!}-\frac{1}{\delta}~
\frac{n+1}{n!}\,s^n\Bigr)\,s'={-}\,\frac{1}{\delta^2}~\frac{s^{n+1}}{n!}~.
\eq
The latter equality simplifies to
\beq \label{e92}
\Bigl(1-\frac{n+1}{\delta}+\frac{s}{\delta}\Bigr)\,s'={-}\,\frac{s}{\delta^2}~,
\eq
and we arrive at (\ref{e29}) where we also use that $s_n(\delta)>n+1-\delta$.\hfill \qedsymbol{} \\
\mbox{}

In the interest of solving the equation in (\ref{e26}), i.e., the equation in (\ref{e87}), we shall show that $\psi''(s)>0$ for $s>n+1-\delta$. For $n=0$, this is obvious since $s-{\rm ln}(1+s/\delta)$ is a strictly convex function of $s>0$. For $n=1,2,...\,$, we compute the numerator $N(s)$ of $\delta\cdot n!\,\psi''(s)$ from (\ref{e88}) as
\newpage
\begin{eqnarray} \label{e93}
& \mbox{} & \hspace*{-8mm}
N(s) \nonumber \\[2mm]
& & \hspace*{-8mm}=~((n+1)\,s^n-n(n+1-\delta)\,s^{n-1})\Bigl(1+ \frac{s^1}{1!}+...+\frac{s^{n-1}}{(n-1)!}+\frac{s^n}{n!}+\frac{s^{n+1}}{\delta\cdot n!}\Bigr) \nonumber \\[3mm]
& & \hspace*{-4mm}-~(s^{n+1}-(n+1-\delta)\,s^n)\Bigl(1+\frac{s^1}{1!}+...+\frac{s^{n-2}}{(n-2)!}+\frac{s^{n-1}} {(n-1)!}+\frac{n+1}{\delta}~\frac{s^n}{n!}\Bigr) \nonumber \\[3mm]
& & \hspace*{-8mm}=~ s^{n-1}((n+1)\,s-n(n+1-\delta))\Bigl(1+\frac{s^1}{1!}+...+\frac{s^{n-1}} {(n-1)!}+\frac{s^n}{n!}+\frac{s^{n+1}}{\delta\cdot n!}\Bigr) \nonumber \\[3mm]
& & \hspace*{-4mm}-~s^n(s-(n+1-\delta))\Bigl(1+\frac{s^1}{1!}+...+\frac{s^{n-2}}{(n-2)!}+ \frac{s^{n-1}}{(n-1)!}+\frac{n+1}{\delta}~\frac{s^n}{n!}\Bigr)~. \nonumber \\[2mm]
\mbox{}
\end{eqnarray}
For $s>n+1-\delta$, we have
\beq \label{e94}
0<n(n+1-\delta)<(n+1)(n+1-\delta)<(n+1)\,s~,
\eq
and so, by the second inequality in (\ref{e94}),
\begin{eqnarray} \label{e95}
& \mbox{} & \hspace*{-8mm}N(s) \nonumber \\[2mm]
& & \hspace*{-8mm}>~s^{n-1}\Bigl[(n+1)(s-(n+1-\delta))\Bigl(1+ \frac{s^1}{1}+...+\frac{s^{n-1}}{(n-1)!}+\frac{s^n}{n!}+\frac{s^{n+1}} {\delta\cdot n!}\Bigr) \nonumber \\[3mm]
& & -~(s-(n+1-\delta))\Bigl(s+\frac{s^2}{1!}+...+\frac{s^{n-1}} {(n-2)!}+\frac{s^n}{(n-1)!}+\frac{n+1}{\delta}~\frac{s^{n+1}}{n!}\Bigr)\Bigr] \nonumber \\[3mm]
& & \hspace*{-8mm}=~s^{n-1}(s-(n+1-\delta))\Bigl[n+1+\frac{n+1-1}{1!}\,s+\frac {n+1-2}{2!}\,s^2 \nonumber \\[3mm]
& & \hspace*{3.8cm}+...+\frac{n+1-n}{(n-1)!}\,s^n+\frac{n+1-(n+1)}{\delta\cdot n!}\,s^{n+1}\Bigr]\,,
\end{eqnarray}
and this is positive. Hence $\psi''(s)>0$ for $s>n+1-\delta$.

The Newton iteration
\beq \label{e96}
s^{(0)}=S\,,~~s^{(j+1)}=s^{(j)}-\frac{\psi(s^{(j)})}{\psi'(s^{(j)})}~,~~~~~~ j=0,1,...
\eq
has iterands $s^{(j)}$ that decrease to $s_n(\delta)$ when the initialization $S$ satisfies $S>s_n(\delta)$. For this, one can choose, following Prop.~\ref{prop8},
\beq \label{e97}
S={\rm ln}\Bigl(\frac{n+1}{\delta}\Bigr)+(n+1-\delta)\,{\rm ln}\Bigl({\rm ln}\Bigl(\frac{n+1}{\delta}\Bigr)\Bigr)+0.8(n+1)~,
\eq
where it has been used that $-0.8\leq\frac{\partial\gamma}{\partial a}\,(1,z)\leq0$, $z\geq0$, see after the proof of Prop.~\ref{prop8} in Sec.~\ref{sec4}.

It is seen from (\ref{e88}) and (\ref{e38}) that the Newton iteration in (\ref{e96}) has problems when $\delta$ is close to $n+1$ and $n$ is not small so that $\psi'(s_n(\delta))$ gets small. For small values of $y=n+1-\delta>0$, one can compute $s_n(\delta)$ by using the B\"urmann-Lagrange inversion theorem. Thus, one writes the equation (\ref{e26}) as
\beq \label{e98}
\frac{1}{n+1-y}=\frac{1}{\delta}=\frac{n!}{s^{n}}\,\Bigl(e^s-1-\frac{s^1}{1!}- ...-\frac{s^n}{n!}\Bigr)=:f(s)~,
\eq
so that
\beq \label{e99}
\frac{(n+1)\,f(s)-1}{f(s)}=y~.
\eq
Now from $\exp(s)=\sum_{k=0}^{\infty}\,s^k/k!$, one gets
\beq \label{e100}
f(s)=\frac{n!}{s^{n+1}}\,\sum_{k=n+1}^{\infty}\,\frac{s^k}{k!}= \sum_{k=0}^{\infty}\,\frac{n!}{(n+k+1)!}\,s^k~,
\eq
\beq \label{e101}
(n+1)\,f(s)-1=\sum_{k=1}^{\infty}\,\frac{(n+1)!}{(n+k+1)!}\,s^k~.
\eq
Some further manipulations bring the equation (\ref{e99}) in the form
\beq \label{e102}
s~\frac{\dsum_{k=0}^{\infty}\,b_k\,s^k}{\dsum_{k=0}^{\infty}\,a_k\,s^k}= \frac{n+2}{n+1}\,y~,
\eq
where
\beq \label{e103}
b_k=\frac{(n+2)!}{(n+k+2)!}\,,~~a_k=\frac{(n+1)!}{(n+k+1)!}~,~~~~~~ k=0,1,...~.
\eq
The B\"urmann-Lagrange formula then gives the small $y$
\beq \label{e104}
s=\sum_{m=1}^{\infty}\,\frac1m\,C_{s^{m-1}}\,\left[\left( \frac{\dsum_{k=0}^{\infty}\,a_k\,s^k}{\dsum_{k=0}^{\infty}\,b_k\,s^k}\right)^m \right]\,\Bigl(\frac{n+2}{n+1}\,y\Bigr)^m~,
\eq
where "$C_{s^{m-1}}[~~]$" is short-hand notation for "the coefficient of $s^{m-1}$ in". To proceed, we write
\beq \label{e105}
\frac{\dsum_{k=0}^{\infty}\,a_k\,s^k}{\dsum_{k=0}^{\infty}\,b_k\,s^k}= \sum_{k=0}^{\infty}\,c_k\,s^k~.
\eq
Since $a_0=1=b_0$, we can compute the $c_k$ recursively according to
\beq \label{e106}
c_0=1\,;~~~~~~c_k=a_k-\sum_{l=0}^{k-1}\,b_{k-l}\,c_l~,~~~~~~k=1,2,...~.
\eq
Having the $c_k$ available, we then compute
\beq \label{e107}
g_m=\frac{1}{m}\,C_{s^{m-1}}\Bigl[\Bigl(\sum_{k=0}^{\infty}\,c_k\,s^k\Bigr)^m\Bigr] ~,~~~~~~m=1,2,...~.
\eq
Thus, $g_1=1$ and for $m=2,3,...$
\begin{eqnarray} \label{e108}
g_m & = & \frac{1}{m}\,C_{s^{m-1}}\,\Bigl[\Bigl(1+\sum_{k=1}^{m-1}\, c_k\,s^k\Bigr)^m\Bigr] \nonumber \\[3mm]
& = & \frac{1}{m}\,C_{s^{m-1}}\,\Bigl[\sum_{j=1}^{m-1}\,\Bigl(\!\ba{c} m \\[1mm] j \ea\!\Bigr)\Bigl(\sum_{k=1}^{m-j}\,c_k\,s^k\Bigr)^j\Bigr]~.
\end{eqnarray}
Then one gets
\begin{eqnarray}
\os{m=2}~~~~~g_2 & = & \tfrac12\,C_s\,\Bigl[\Bigl(\!\ba{c} 2 \\[1mm] 1 \ea\!\Bigr)\, c_1\,s\Bigr]=c_1~, \nonumber \\[3mm]
\os{m=3}~~~~~g_3 & = & \tfrac13\,C_{s^2}\,\Bigl[\Bigl(\!\ba{c} 3 \\[1mm] 1 \ea\!\Bigr)\,(c_1\,s+c_2\,s^2)+\tfrac32\,(c_1\,s)^2\Bigr] \nonumber \\[3mm]
& = & \tfrac13\,(3c_2+3c_1^2)=c_2+c_1^2~, \nonumber
\end{eqnarray}
and in a similar fashion
\begin{eqnarray}
\os{m=4}~~~~~g_4 & = & c_3+3c_1c_2+c_1^3~, \nonumber \\[3mm]
\os{m=5}~~~~~g_5 & = & c_4+4c_1c_3+2c_2^2+6c_1^2c_2+c_1^4~. \nonumber
\end{eqnarray}
For example, with $n=1$ and $y=2-\delta$, we find
\beq \label{e109}
s_1(\delta)=\tfrac32\,y+\tfrac{1}{12}\,(\tfrac32\,y)^2+\tfrac{7}{360}\,(\tfrac32\,y)^3 +\tfrac{41}{8640}\,(\tfrac32\,y)^4+\tfrac{2243}{1814400}\,(\tfrac32\,y)^5+...~,
\eq
and, with $n=3$ and $y=4-\delta$,
\beq \label{e110}
s_3(\delta)=\tfrac54\,y+\tfrac{1}{30}\,(\tfrac54\,y)^2+\tfrac{8}{1575}\,(\tfrac54\,y)^3+\tfrac{289}{378000} \,(\tfrac54\,y)^4+\tfrac{1181}{9922500}\,(\tfrac54\,y)^5+...~.
\eq

\section{Bounds for $s_n$} \label{sec4}
\mbox{} \\[-9mm]

According to Taylor's formula with remainder in integral form, see \cite{ref3}, 1.4(vi), we have for $f\in C^{n+1}[0,\infty)$
\beq \label{e111}
f(s)=\sum_{k=0}^n\,\frac{f^{(k)}(0)}{k!}\,s^k+\frac{1}{n!}\,\il_0^s\,(s-t)^n\, f^{(n+1)}(t)\,dt~,~~~~~~s>0~.
\eq
Using this with $f(s)=\exp(s)$ and substituting $x=s-t\in[0,s]$ in the integral at the right-hand side of (\ref{e111}), we get
\beq \label{e112}
e^s=1+\frac{s^1}{1!}+...+\frac{s^n}{n!}+\frac{1}{n!}\,e^s\,\il_0^s\,x^n\,e^{-x}\,dx ~,~~~~~~s>0~.
\eq
Combining this with Prop.~\ref{prop4}, we get for $n=0,1,...\,$, $\delta\in(0,n+1)$ and $s>0$
\beq \label{e113}
s\:>~{\rm or}~<\:s_n(\delta)\Leftrightarrow\il_0^s\,x^n\, e^{-x}\,dx\:>~{\rm or}~<\:\frac{1}{\delta}\,s^{n+1}\,e^{-s}~.
\eq
\mbox{} \\
{\bf Proof of Proposition~\ref{prop5}.} \mbox{}

a.~~Let $n=0,1,...\,$, and assume that $S\in C^1(0,n+1]$ with $S(n+1)=0<S(\delta)$, $\delta\in(0,n+1)$. According to (\ref{e113}), we have $S(\delta)<s_n(\delta)$ if and only if
\beq \label{e114}
\il_0^{S(\delta)}\,x^n\,e^{-x}\,dx<\frac{1}{\delta}\,(S(\delta))^{n+1}\, e^{-S(\delta)}~.
\eq
We have equality in (\ref{e114}) when $\delta=n+1$ since $S(n+1)=0$. Therefore (\ref{e114}) holds for $\delta\in(0,n+1)$ when
\beq \label{e115}
\frac{d}{d\delta}\,\left[\il_0^{S(\delta)}\,x^n\,e^{-x}\,dx\right]> \frac{d}{d\delta}\,\Bigl[\frac{1}{\delta}\,(S(\delta))^{n+1}\,e^{-S(\delta)} \Bigr]~,~~~~~~\delta\in(0,n+1)~.
\eq
A simple calculation shows that (\ref{e115}) holds if and only if
\beq \label{e116}
S'(\delta)>{-}\,\frac{1}{\delta^2}\,S(\delta)+\frac{1}{\delta}\,S'(\delta)(n+1-S(\delta)) ~,~~~~~~\delta\in(0,n+1)~,
\eq
where we have used that $S(\delta)>0$, $\delta\in(0,n+1)$. The condition (\ref{e116}) is equivalent with (\ref{e35}).

b.~~Now assume that $S\in C^1(0,n+1]$ with $S(n+1)\geq0$, $S(\delta)>0$, $\delta\in(0,n+1)$. Also assuming (\ref{e36}), we now want to show that
\beq \label{e117}
\il_0^{S(\delta)}\,x^n\,e^{-x}\,dx>\frac{1}{\delta}\,(S(\delta))^{n+1}\,e^{-S(\delta)}~,~~~~~~ \delta\in(0,n+1)~.
\eq
In the case that $S(n+1)=0$, we can argue in exactly the same way as in the proof of a, with reversed inequality signs in (\ref{e115}) and (\ref{e116}). In the case that $S(n+1)>0$, we have that (\ref{e117}) also holds for $\delta=n+1$, since
\beq \label{e118}
\il_0^S\,x^n\,e^{-x}\,dx>e^{-S}\,\il_0^S\,x^n\,dx=\frac{1}{n+1}\, S^{n+1}\,e^{-S}~,~~~~~~S>0~,
\eq
and so we are done as well.\hfill \qedsymbol{} \\ \\
{\bf Proof of Proposition~\ref{prop6}.}~~To show the second inequality in (\ref{e37}), we take
\beq \label{e119}
S(\delta)=(n+1-\delta)+(n+1-\delta)/\delta~,~~~~~~\delta\in(0,n+1]~.
\eq
Observe that $S(n+1)=0$. When checking (\ref{e36}), it is noted that a factor $(n+1-\delta)$ cancels, and it remains to check that
\beq \label{e120}
1+\frac{1}{\delta}<{-}S'(\delta)~,~~~~~~\delta\in(0,n+1)~.
\eq
Since $S'(\delta)={-}1-(n+1)/\delta^2$, we see that (\ref{e120}) obviously holds.

To show the first inequality in (\ref{e37}), we take
\beq \label{e121}
S(\delta)={\rm ln}\Bigl(\frac{n+1}{\delta}\Bigr)+(n+1-\delta)~,~~~~~~ \delta\in(0,n+1]~.
\eq
Observe that $S(n+1)=0$. Since $S'(\delta)={-}1/\delta-1$, the condition (\ref{e35}) is satisfied if and only if
\beq \label{e122}
{\rm ln}\Bigl(\frac{n+1}{\delta}\Bigr)+(n+1-\delta)>(\delta+1)\, {\rm ln}\Bigl(\frac{n+1}{\delta}\Bigr)~,~~~~~~\delta\in(0,n+1)~,
\eq
i.e., if and only if
\beq \label{e123}
{\rm ln}\Bigl(\frac{n+1}{\delta}\Bigr)<\frac{n+1}{\delta}-1~,~~~~~~\delta\in (0,n+1)~.
\eq
Since $(n+1)/\delta>1$ for $\delta\in(0,n+1)$, we see that (\ref{e123}) indeed holds.\hfill \qedsymbol{} \\
\mbox{}

The next result is required in the proof of convexity of $s_n(\delta)$ as a function of $\delta\in(0,n+1)$.

\setcounter{thm}{0}
\begin{lem} \label{lem1}
Let $n=0,1,...\,$. Then
\beq \label{e124}
s_n(\delta)>(n+1-\delta)+((n+1-\delta)+\tfrac14\,\delta^2)^{1/2}-\tfrac12\, \delta~,~~~~~~\delta\in(0,n+1)~.
\eq
\end{lem}

\noindent
{\bf Proof.}~~Denoting the right-hand side of (\ref{e124}) by $S(\delta)$, we have with $y=n+1-\delta$
\beq \label{e125}
S(\delta)=y+(y+\tfrac14\,\delta^2)^{1/2}-\tfrac12\,\delta~,~~~~~~ \delta\in(0,n+1)~.
\eq
We have $S(n+1)=0$, and we shall show that (\ref{e35}) holds, i.e., that
\beq \label{e126}
\frac{1}{\delta}\,S(\delta)>{-}S'(\delta)(S(\delta)-y)~,~~~~~~ \delta\in(0,n+1)~.
\eq
Inserting
\beq \label{e127}
S(\delta)=y+\frac{y}{(y+\tfrac14\,\delta^2)^{1/2}+\tfrac12\,\delta}
\eq
into (\ref{e126}), we see that we must show that
\beq \label{e128}
-\delta\,S'(\delta)<(y+\tfrac14\,\delta^2)^{1/2}+\tfrac12\,\delta+1~,~~~~~~ \delta\in(0,n+1)~.
\eq
Here we have omitted a few simple manipulations. Since $-(\delta\, S(\delta))'= {-}\delta\,S'(\delta)-S(\delta)$, we see that (\ref{e128}) is equivalent with
\beq \label{e129}
-(\delta\,S'(\delta))'<(y+\tfrac14\,\delta^2)^{1/2}+\tfrac12\,\delta+1-S(\delta)~ ,~~~~~~\delta\in(0,n+1)~.
\eq
Using (\ref{e125}), we see that the right-hand side of (\ref{e129}) equals $\delta-y+1={-}n+2\delta$. Next, using (\ref{e125}) one more time, we get
\begin{eqnarray} \label{e130}
(\delta\,S(\delta))' & = & (y\delta+\delta(y+\tfrac14\,\delta^2)^{1/2}-\tfrac12\,\delta^2)' \nonumber \\[3mm]
& = & n+1-3\delta+(y+\tfrac14\,\delta^2)^{1/2}+\tfrac12\,\delta(\tfrac12\,\delta -1)(y+\tfrac14\,\delta^2)^{-1/2}~, \nonumber \\
\mbox{}
\end{eqnarray}
where we have also used that $y=n+1-\delta$, $y'=1$. The right-hand side of (\ref{e129}) equals $2\delta-n$, and so the condition in (\ref{e129}) is equivalent with
\beq \label{e131}
n+1-3\delta+(y+\tfrac14\,\delta^2)^{1/2}+\tfrac12\,\delta(\tfrac12\,\delta-1) (y+\tfrac14\,\delta^2)^{-1/2}>n-2\delta~,~~~~~~ \delta\in(0,n+1)~,
\eq
i.e., with
\beq \label{e132}
(y+\tfrac14\,\delta^2)^{1/2}+\tfrac12\,\delta(\tfrac12\,\delta-1)(y+\tfrac14\,\delta^2)^{-1/2} >\delta-1~,~~~~~~\delta\in(0,n+1)~,
\eq
i.e., with
\beq \label{e133}
y-\tfrac12\,\delta+\tfrac12\,\delta^2>(\delta-1)(y+\tfrac14\,\delta^2)^{1/2}~,~~~~~~ \delta\in(0,n+1)~.
\eq
When $0<\delta<1$, the right-hand side of (\ref{e133}) is negative while the left-hand side of (\ref{e133}) equals $n+\frac12\,(1-\delta)(2-\delta)>0$, and so (\ref{e133}) holds. Therefore, assume $\delta\geq1$, which implies that $n\geq1$ since $\delta\in(0,n+1)$. We thus need to check whether
\beq \label{e134}
(y-\tfrac12\,\delta+\tfrac12\,\delta^2)^2>(\delta-1)^2(y+\tfrac14\,\delta^2)~,
\eq
i.e., upon expanding either side of (\ref{e134}) whether $y+\delta>1$. Since $y=n+1-\delta$ and $n\geq1$, we have $y+\delta>1$, indeed.\hfill \qedsymbol{} \\ \\
{\bf Proof of Proposition~\ref{prop7}.}~~We must check that $s_n''(\delta)>0$, $\delta\in(0,n+1)$, where $n=0,1,...\,$. We have by Prop.~\ref{prop3} for $\delta\in(0,n+1)$
\begin{eqnarray} \label{e135}
s_n''(\delta) & = & \frac{d}{d\delta}\,\Bigl[{-}\,
\frac{1}{\delta}~\frac{s_n(\delta)} {s_n(\delta)-(n+1-\delta)}\Bigr] \nonumber \\[3mm]
& = & \frac{1}{\delta^2}~\frac{s_n(\delta)}{s_n(\delta)-(n+1-\delta)} +\frac{1}{\delta}~\frac{s_n(\delta)}{(s_n(\delta)-(n+1-\delta))^2} \nonumber \\[3mm]
& & \hspace*{4.2cm} +~\frac{n+1-\delta}{\delta}~\frac{s_n'(\delta)}{(s_n(\delta)-(n+1-\delta))^2} \nonumber \\[3mm]
& = & \frac{1}{\delta^2}~\frac{s_n(\delta)} {s_n(\delta)-(n+1-\delta)}\, \Bigl(1+\frac{\delta}{s_n(\delta)-(n+1-\delta)} \nonumber \\[3mm]
& & \hspace*{4.4cm}-~\frac{n+1-\delta}{(s_n(\delta)-(n+1-\delta))^2}\Bigr)~,
\end{eqnarray}
where in the last step Prop.~\ref{prop3} has been used once more. Thus, with $y=n+1-\delta$, we have $s_n''(\delta)>0$ if and only if
\beq \label{e136}
1+\frac{\delta}{s_n(\delta)-y}-\frac{y}{(s_n(\delta)-y)^2}>0~.
\eq
A simple computation shows that (\ref{e136}) is equivalent with
\beq \label{e137}
(s_n(\delta)-y+\tfrac12\,\delta)^2>y+\tfrac14\,\delta^2~.
\eq
Then Lemma~\ref{lem1} gives the result.\hfill\qedsymbol{} \\
\mbox{}

The upper bound in Prop.~\ref{prop6} can be sharpened as follows.

\begin{lem} \label{lem2}
Let $n=0,1,...\,$, and $A>1$, $B=A/(A-1)$. Then
\beq \label{e138}
s_n(\delta)<A\,{\rm ln}\Bigl(\frac{n+1}{\delta}\Bigr) +B(n+1-\delta)~,~~~~~~ \delta\in(0,n+1)~.
\eq
\end{lem}

\noindent
{\bf Proof.}~~Denote the right-hand side of (\ref{e138}) by $S(\delta)$. Then $S(n+1)=0$, and so it suffices by Prop.~\ref{prop5}b to check whether
\beq \label{e139}
S(\delta)<{-}\delta\,S'(\delta)(S(\delta)-(n+1-\delta))~,~~~~~~ \delta\in(0,n+1)~.
\eq
Since $S'(\delta)={-}A/\delta-B$, we have that (\ref{e139}) holds if and only if
\beq \label{e140}
(B-(A+B\delta)(B-1))(n+1-\delta)<(A-1+B\delta)\,{\rm ln}\Bigl(\frac{n+1}{\delta} \Bigr)\,,~~~~\delta\in(0,n+1)~.
\eq
With $A>1$ and $B=A/(A-1)$, we have
\beq \label{e141}
B-(A+B\delta)(B-1)={-}B(B-1)\,\delta<0\,,~~ A-1+B\delta>0\,,~~~~\delta\in (0,n+1)~,
\eq
and so (\ref{e140}) is valid.\hfill\qedsymbol{} \\ \\
When we choose $A=2$ in Lemma~\ref{lem2}, so that $B=2$ as well, the right-hand side of (\ref{e138}) equals twice the lower bound for $s_n(\delta)$ in the left-hand side of (\ref{e37}) in Prop.\ \ref{prop6}. From the first inequality  in (\ref{e37}) and Lemma~\ref{lem2} we have
\beq \label{e142}
\lim_{\delta\downarrow0}\:\frac{s_n(\delta)}{{\rm ln}\Bigl(\dfrac{n+1}{\delta}\Bigr)}=1
\eq
for $n=0,1,...\,$. We aim at further sharpening of (\ref{e142}) in which we are guided by the asymptotics of $s_0(\delta)$ as $s\downarrow0$.

\begin{lem} \label{lem3}
We have
\beq \label{e143}
s_0(\delta)={\rm ln}\Bigl(\frac{1}{\delta}\Bigr)+{\rm ln}\Bigl({\rm ln}\Bigl( \frac{1}{\delta}\Bigr)\Bigr)+O\left(\frac{{\rm ln}\Bigl({\rm ln}\Bigl(\dfrac{1}{\delta}\Bigr)\Bigr)}{{\rm ln}\Bigl(\dfrac{1}{\delta}\Bigr)} \right) ~,~~~~~~\delta\downarrow0~.
\eq
\end{lem}

\noindent
{\bf Proof.}~~By Prop.~\ref{prop1} we have that $s_0(\delta)$ is the unique solution $s>0$ of the equation
\beq \label{e144}
e^s=1+\frac{s}{\delta}~,~~~~~~\delta\in(0,1)~.
\eq
Furthermore, we have from (\ref{e142}) that $s_0(\delta)={\rm ln}(\frac{1}{\delta})(1+o(1))$, $\delta\downarrow0$. From (\ref{e144}) we get for $\delta\downarrow0$
\begin{eqnarray} \label{e145}
s & = & {\rm ln}\Bigl(1+\frac{s}{\delta}\Bigr)={\rm ln}\Bigl(1+\frac{1}{\delta}\,{\rm ln}\Bigl(1+\frac{s}{\delta}\Bigr)\Bigr) \nonumber \\[3mm]
& = & {\rm ln}\Bigl(1+\frac{1}{\delta}\,{\rm ln}\Bigl(\frac{s}{\delta}\Bigr)+\frac{1}{\delta}\,{\rm ln}\Bigl(1+\frac {\delta}{s}\Bigr)\Bigr) \nonumber \\[3mm]
& = & {\rm ln}\Bigl(1+\frac{1}{\delta}\,{\rm ln}\Bigl(\frac{s}{\delta}\Bigr)+O\Bigl(\frac1s\Bigr)\Bigr)
={\rm ln}\Bigl(\frac{1}{\delta}\,{\rm ln}\Bigl(\frac{s}{\delta}\Bigr)\Bigr) +O\left(\frac{\delta}{s\,{\rm ln}\Bigl(\dfrac{s}{\delta}\Bigr)}\right) \nonumber \\[3mm]
& = & {\rm ln}\Bigl(\frac{1}{\delta}\Bigr)+{\rm ln}\Bigl({\rm ln}\,s+{\rm ln}\Bigl(\frac{1}{\delta}\Bigr)\Bigr)+O\left(\frac{\delta}{s\,{\rm ln}\Bigl(\dfrac{s}{\delta}\Bigr)}\right) \nonumber \\[3mm]
& = & {\rm ln}\Bigl(\frac{1}{\delta}\Bigr)+{\rm ln}\Bigl({\rm ln}\Bigl(\frac{1}{\delta}\Bigr)\Bigr)+O\left(\frac{{\rm ln}\,s}{{\rm ln}\Bigl(\dfrac{1}{\delta}\Bigr)}\right)+O\left(\frac{\delta}{s\,{\rm ln}\,\Bigl(\dfrac{s}{\delta}\Bigr)}\right)~,
\end{eqnarray}
and the result follows from $s={\rm ln}(\frac{1}{\delta})(1+o(1))$, $\delta\downarrow0$.\hfill\qedsymbol{} \\ \\
A result like Lemma~\ref{lem3} can be given for $s_n(\delta)$, $n=1,2,...\,$, where now a leading behaviour
\beq \label{e146}
{\rm ln}\Bigl(\frac{n+1}{\delta}\Bigr)+(n+1)\,{\rm ln}\Bigl({\rm ln}\Bigl( \frac{n+1}{\delta}\Bigr)\Bigr)
\eq
appears as $\delta\downarrow0$. Note that we need to consider $\delta\downarrow0$ in Lemma~\ref{lem3} and (\ref{e146}) while we aim at bounds that are valid for all $\delta\in(0,n+1)$. \\ \\
{\bf Proof of Proposition~\ref{prop8}.}~~We have from Prop.~\ref{prop3} and the first inequality in (\ref{e37})
\beq \label{e147}
s_n'(\delta)={-}\frac{1}{\delta}-\frac{n+1-\delta}{s_n(\delta)-(n+1-\delta)}> {-}\,\frac{1}{\delta}-\frac{1}{\delta}~\frac{n+1-\delta}{{\rm ln}\Bigl(\dfrac{n+1}{\delta}\Bigr)}~,~~~~~~\delta\in(0,n+1)~.
\eq
From $s_n(n+1)=0$, we then get
\begin{eqnarray} \label{e148}
s_n(\delta) & = & {-}(s_n(n+1)-s_n(\delta))={-}\,\il_{\delta}^{n+1}\, s_n'(\delta_1)\,d\delta_1 \nonumber \\[3mm]
& < & \il_{\delta}^{n+1}\,\Bigl(\frac{1}{\delta_1}+\frac{1}{\delta_1}~ \frac{n+1-\delta_1}{{\rm ln}\Bigl(\dfrac{n+1}{\delta_1}\Bigr)}\Bigr)\, d\delta_1 \nonumber \\[3mm]
& = & {\rm ln}\Bigl(\frac{n+1}{\delta}\Bigr)+\il_{\delta}^{n+1}\,\frac{1}{\delta_1} ~\frac{n+1-\delta_1}{{\rm ln}\Bigl(\dfrac{n+1}{\delta_1}\Bigr)}\,d\delta_1 ~,~~~~~~\delta\in(0,n+1)~. \nonumber \\[-4mm]
\mbox{}
\end{eqnarray}
We have for the remaining integral
\beq \label{e149}
\il_{\delta}^{n+1}\,\frac{1}{\delta_1}~\frac{n+1-\delta_1} {{\rm ln}\Bigl(\dfrac{n+1}{\delta_1}\Bigr)}\,d\delta_1=(n+1)\,\il_1^x\, \frac{x_1-1}{x_1^2\,{\rm ln}\,x_1}\,dx_1~,~~~~~~x=\frac{n+1}{\delta}~,
\eq
where we have substituted $x_1=(n+1)/\delta_1\in(1,x)$. Next, it is straightforward to check by differentiation that for $x>1$
\beq \label{e150}
\il_1^x\,\frac{x_1-1}{x_1^2\,{\rm ln}\,x_1}\,dx_1=\Bigl(1-\frac1x\Bigr)\,{\rm ln}({\rm ln}\,x)-\il_0^{{\rm ln}\,x}\,e^{-t}\,{\rm ln}\,t\,dt~.
\eq
Hence, going back to (\ref{e149}) and (\ref{e148}), we get for $\delta\in(0,n+1)$
\beq \label{e151}
s_n(\delta)<z+(n+1-\delta)\,{\rm ln}\,z-(n+1)\,\il_0^z\,e^{-t}\,{\rm ln}\,t\,dt~,
\eq
where
\beq \label{e152}
z={\rm ln}\,x={\rm ln}\Bigl(\frac{n+1}{\delta}\Bigr)~.
\eq
For the remaining integral in (\ref{e151}), we have
\beq \label{e153}
\Phi(z):=\il_0^z\,e^{-t}\,{\rm ln}\,t\,dt=\frac{d}{da}\,\gamma(a,z)\Bigr|_{a=1}~, ~~~~~~z\geq0~,
\eq
where $\gamma(a,z)$ is the incomplete $\Gamma$-function, \cite{ref3}, Ch.~8,
\beq \label{e154}
\gamma(a,z)=\il_0^z\,t^{a-1}\,e^{-t}\,dt~,~~~~~~a>0\,,~~z>0~.
\eq
This gives (\ref{e40}--\ref{e41}).\hfill\qedsymbol{} \\ \\
Observe that the two functions $(n+1-\delta)\,{\rm ln}\,z$ and $\frac{\partial\gamma}{\partial a}\,(1,z)$ with $z={\rm ln}((n+1)/\delta)$ are slightly singular as $\delta\uparrow(n+1)$ while the right-hand side of (\ref{e40}) is regular as $\delta\uparrow(n+1)$. With $y=n+1-\delta$, we have for the right-hand side of (\ref{e40}) the expansion, compare (\ref{e38}--\ref{e39})
\beq \label{e155}
\frac{n+2}{n+1}\,y+\frac{n+3}{4(n+1)^2}\,y^2+\frac{5n+17}{36(n+1)^3}\,y^3+...~.
\eq
The integral at the left-hand side of (\ref{e150}) is an analytic function of $x$, $|x-1|<1$; for the sake of identification in terms of (special) functions we have written this analytic function in (\ref{e150}) as a difference of two functions that both exhibit a $(x-1)\,{\rm ln}(x-1)$-behaviour as $x\downarrow1$.

We present some more information on the function $\Phi(z)$ in (\ref{e153}). We have for $z\geq0$
\beq \label{e156}
0=\Phi(0)>\Phi(z)\geq\Phi(1)={-}\,\sum_{k=0}^{\infty}\,\frac{({-}1)^k}{k!}~ \frac{1}{(k+1)^2}={-}0.796599599...\,,
\eq
see \cite{ref4}, A001563 and \cite{ref5}, eq.~289. Furthermore, $\Phi(z)\pr {-}\gamma={-}0.5772156649...\,$, where $\gamma$ is Euler's gamma. Finally, by partial integration in (\ref{e153}), we have for $z>0$
\beq \label{e157}
\Phi(z)=(1-e^{-z})\,{\rm ln}\,z-\il_0^z\,\frac{1-e^{-t}}{t}\,dt= (1-e^{-z})\,{\rm ln}\,z-{\rm Ein}(z)~,
\eq
where ${\rm Ein}(z)$ is the complementary exponential integral, see \cite{ref3}, Ch.~6. Consequently, by \cite{ref3}, 6.2.3 and 6.6.4, we get
\beq \label{e158}
\Phi(z)=(1-e^{-z})\,{\rm ln}\,z+\sum_{k=1}^{\infty}\, \frac{({-}z)^k}{k\cdot k!}~,~~~~~~z\geq0~.
\eq

\newpage
\noindent
\section{The function $ME_{n,\delta}$} \label{sec5}
\mbox{} \\[-9mm]

We next consider the function $ME_{n,\delta}$, $\delta\in(0,n+1)$. \\ \\
{\bf Proof of Proposition~\ref{prop9}.}~~To show that $ME_{n,\delta}\pr1$ as $\delta\downarrow0$, we observe from Props.~\ref{prop6} and \ref{prop8} that
\beq \label{e159}
s_n(\delta)={\rm ln}\Bigl(
\frac{n+1}{\delta}\Bigr)+O\Bigl({\rm ln}\Bigl({\rm ln}\Bigl(\frac{n+1}{\delta}\Bigr)\Bigr)\Bigr)~,~~~~~~\delta\downarrow0~.
\eq
Hence,
\beq \label{e160}
(s_n(\delta))^{\delta}=1+o(1)~,~~~~~~ \Bigl(1+\frac{s^1}{1!}+...+\frac{s^n}{n!}\Bigr)\,e^{-s}\Bigr|_{s=s_n(\delta)} =o(1)
\eq
as $\delta\downarrow0$. From the first expression for $ME_{n,\delta}$ in (\ref{e28}) of Prop.~\ref{prop2} we then get that $ME_{n,\delta}\pr1$, $\delta\downarrow0$.

Next, from (\ref{e37}--\ref{e39}), we have
\beq \label{e161}
s_n(\delta)=\frac{n+2}{n+1}\,(n+1-\delta)+O((n+1-\delta)^2)~,~~~~~~ \delta\uparrow(n+1)~.
\eq
Hence,
\beq \label{e162}
(s_n(\delta))^{n+1-\delta}=1+o(1)~,~~~~~~ e^{s_n(\delta)}=1+o(1)
\eq
as $\delta\uparrow(n+1)$. From the second expression for $ME_{n,\delta}$ in (\ref{e28}) of Prop.~\ref{prop2} we then get that $ME_{n,\delta}\pr 1/(n+1)!$, $\delta\uparrow(n+1)$.\hfill\qedsymbol{} \\ \\
{\bf Proof of Proposition~\ref{prop10}.}~~We have by Prop.~\ref{prop2}
\beq \label{e163}
ME_{n,\delta}=\frac{s^{n+1-\delta}}{n!\,\delta\,e^s}~,~~~~~~ s=s_n(\delta)\,,~~\delta\in(0,n+1)~.
\eq
Hence, by Prop.~\ref{prop3},
\begin{eqnarray} \label{e164}
& \mbox{} & \frac{d}{d\delta}\,[{\rm ln}\,ME_{n,\delta}] = \frac{d}{d\delta}\,[(n+1-\delta)\,{\rm ln}(s_n(\delta))-{\rm ln}\,\delta-s_n(\delta)] \nonumber \\[3mm]
& & =~{-}{\rm ln}(s_n(\delta))+(n+1-\delta)\,\frac{s_n'(\delta)}{s_n(\delta)} -\frac{1}{\delta}-s_n'(\delta) \nonumber \\[3mm]
& & =~{-}{\rm ln}(s_n(\delta))-\frac{1}{\delta}~\frac{n+1-\delta} {s_n(\delta)-(n+1-\delta)}-\frac{1}{\delta}+\frac{1}{\delta}~ \frac{s_n(\delta)}{s_n(\delta)-(n+1-\delta)} \nonumber \\[3mm]
& & =~{-}{\rm ln}(s_n(\delta))~,
\end{eqnarray}
and this is (\ref{e45}).\hfill\qedsymbol{} \\ \\
{\bf Proof of Proposition~\ref{prop11}.}~~From Prop.~\ref{prop10} we have that the unique minimum of the continuous and strictly convex function $ME_{n,\delta}$, $\delta\in[0,n+1]$, is assumed at the unique $\delta\in(0,n+1)$ such that $s_n(\delta)=1$. Here it has also been used that the continuous function $s_n(\delta)$ decreases strictly from $+\infty$ at $\delta=0$ to 0 at $\delta=n+1$ Denote this unique $\delta$ by $\delta_{n,E}$ and the minimum value of $ME_{n,\delta}$ by $\hat{E}_n$. We have from
\beq \label{e165}
e^s=1+\frac{s^1}{1!}+...+\frac{s^n}{n!}+\frac{1}{\delta}~\frac{s^{n+1}}{n!}~, ~~~~~~s=s_n(\delta)\,,~~\delta\in(0,n+1)~,
\eq
see Prop.~\ref{prop1}, and $s_n(\delta_{n,E})=1$ that
\beq \label{e166}
e=1+\frac{1}{1!}+...+\frac{1}{n!}+\frac{1}{\delta_{n,E}}~\frac{1}{n!}~.
\eq
Therefore,
\beq \label{e167}
\delta_{n,E}=\Bigl(n!\Bigl(e-\sum_{k=0}^n\,\frac{1}{k!}\Bigr)\Bigr)^{-1}~.
\eq
Furthermore, we have from
\beq \label{e168}
ME_{n,\delta}=\frac{s^{n+1-\delta}}{n!\,\delta\,e^s}~,~~~~~~ s=s_n(\delta)\,,~~\delta\in(0,n+1)~,
\eq
see Prop.~\ref{prop2}, and $s_n(\delta_{n,E})=1$ that
\beq \label{e169}
\hat{E}_n=ME_{n,\delta_{n,E}}=\frac{1}{n!\,\delta_{n,E}\,e}=\frac1e\,\Bigl( e-\sum_{k=0}^n\,\frac{1}{k!}\Bigr)~.
\eq
Thus, we have (\ref{e46}) and (\ref{e47}) as required.\hfill\qedsymbol{} \\
\mbox{}

We shall now give bounds for $\delta_{n,E}$ and $\hat{E}_n$.

\begin{lem} \label{lem4}
Let $n=0,1,...\,$. Then
\begin{eqnarray} \label{e170}
& \mbox{} & \frac{n+1}{n+2}\,(n+1)<\delta_{n,E}<\frac{n+2}{n+3}\,(n+1)~, \nonumber \\[3mm]
& & \frac{n+3}{n+2}~\frac{1}{(n+1)!\,e}<\hat{E}_n<\frac{n+2}{n+1}~ \frac{1}{(n+1)!\,e}~.
\end{eqnarray}
\end{lem}

\noindent
{\bf Proof.}~~We have
\begin{eqnarray} \label{e171}
& \mbox{} & n!\Bigl(e-\sum_{k=0}^n\,\frac{1}{k!}\Bigr)=\sum_{k=n+1}^{\infty}\, \frac{n!}{k!} \nonumber \\[3mm]
& & =~\frac{1}{n+1}+\frac{1}{(n+1)(n+2)}+\frac{1}{(n+1)(n+2)(n+3)}+...~.
\end{eqnarray}
Hence,
\beq \label{e172}
\sum_{k=n+1}^{\infty}\,\frac{n!}{k!}>\frac{1}{n+1}\,\Bigl(1+\frac{1}{n+2}\Bigr) =\frac{1}{n+1}~\frac{n+3}{n+2}~,
\eq
and
\begin{eqnarray} \label{e173}
\sum_{k=n+1}^{\infty}\,\frac{n!}{k!} & = & \frac{1}{n+1}\, \Bigl(1+\frac{1}{n+2}+\frac{1}{(n+2)(n+3)}+...\Bigr) \nonumber \\[3mm]
& < & \frac{1}{n+1}\,\Bigl(1+\frac{1}{n+2}+\frac{1}{(n+2)^2}+...\Bigr) \nonumber \\[3mm]
& = & \frac{1}{n+1}~\frac{1}{1-\dfrac{1}{n+2}}=\frac{1}{n+1}~\frac{n+2}{n+1}~.
\end{eqnarray}
Then (\ref{e170}) follows from (\ref{e46}) and (\ref{e47}).\hfill\qedsymbol{}

\section{Lower and upper bounds for $ME_{n,\delta}$} \label{sec6}
\mbox{} \\[-9mm]

At the end of Subsec.~\ref{subsec2.1} we have given a method to pass from a lower bound $S(\delta)$ of $s_n(\delta)$ with $S(\delta)>n+1-\delta$ to two upper bounds for $ME_{n,\delta}$, using the functions $F_1$ and $F_2$ given in (\ref{e48}), for which Prop.~\ref{prop12} is crucial. \\ \\
{\bf Proof of Proposition~\ref{prop12}.}~~Let $n=0,1,...\,$, and let $\delta\in(0,n+1)$. We must show that both $F_1$ and $F_2$ are strictly decreasing in $s\in(n+1-\delta,s_n(\delta))$. This is obvious for $F_1$ since
\beq \label{e174}
\frac{d}{ds}\,[e^{-s}\,s^{n+1-\delta}]=e^{-s}\,s^{n-\delta}(n+1-\delta-s)~,
\eq
so that $F_1'(s)>0$, $s\in(0,n+1-\delta)$ and $F_1'(s)<0$, $s\in(n+1-\delta,\infty)$.

As to $F_2$, we consider for $s>0$
\beq \label{e175}
M(s):={\rm ln}[\delta\cdot n!\,F_2(s)] =(n+1-\delta)\,{\rm ln}\,s-{\rm ln}\, \Bigl[1+\frac{s^1}{1!}+...+\frac{s^n}{n!}+\frac{1}{\delta}~\frac{s^{n+1}}{n!} \Bigr]~.
\eq
We have from a computation as in Prop.~\ref{prop1}, see (\ref{e87}--\ref{e88}), for $s>0$
\begin{eqnarray} \label{e176}
M'(s) & = & \frac{n+1-\delta}{s}-\frac{1+\dfrac{s^1}{1!}+...+\dfrac{s^{n-1}} {(n-1)!}+\dfrac{n+1}{\delta}~\dfrac{s^n}{n!}} {1+\dfrac{s^1}{1!}+...+ \dfrac{s^n}{n!}+\dfrac{1}{\delta}~\dfrac{s^{n+1}}{n!}} \nonumber \\[4mm]
& = & \frac{n+1-\delta}{s}-1+1-\frac{1+\dfrac{s^1}{1!}+...+\dfrac{s^{n-1}}{(n-1)!} +\dfrac{n+1}{\delta}~\dfrac{s^n}{n!}} {1+\dfrac{s^1}{1!}+...+\dfrac{s^n}{n!}+\dfrac{1}{\delta}~\dfrac{s^{n+1}}{n!}} \nonumber \\[4mm]
& = & \frac{(n+1-\delta)-s}{s}+\frac{\dfrac{s^n}{n!}+\dfrac{1}{\delta}~ \dfrac{s^{n+1}}{n!}-\dfrac{n+1}{\delta}~\dfrac{s^n}{n!}} {1+\dfrac{s^1}{1!}+...+\dfrac{s^n}{n!}+\dfrac{1}{\delta}~\dfrac{s^{n+1}}{n!}} \nonumber \\[4mm]
& = & \frac{(n+1-\delta)-s}{s}+\frac{s^n}{\delta\cdot n!}~\frac {s-(n+1-\delta)} {1+\dfrac{s^1}{1!}+...+\dfrac{s^n}{n!}+\dfrac{1}{\delta}~\dfrac{s^{n+1}}{n!}} \nonumber \\[4mm]
& = & \frac{(n+1-\delta)-s}{s}\,\left( 1-\frac{\dfrac{s^{n+1}} {\delta\cdot n!}} {1+\dfrac{s^1}{1!}+...+\dfrac{s^n}{n!}+\dfrac{s^{n+1}}{\delta\cdot n!}} \right)~.
\end{eqnarray}
The factor between (~~) on the last line of (\ref{e176}) is strictly between 0 and 1 for positive $s$. Hence, $M'(s)>0$, $s\in(0,n+1-\delta)$ and $M'(s)<0$ for $s\in(n+1-\delta,\infty)$.

We finally show that $F_2(s)<F_1(s)$, $s\in(n+1-\delta,s_n(\delta))$. We have for $s>0$ that $F_2(s)<F_1(s)$ if and only if
\beq \label{177}
1+\frac{s^1}{1!}+...+\frac{s^n}{n!}+\frac{1}{\delta}~\frac{s^{n+1}}{n!}>e^s~,
\eq
see (\ref{e48}). Hence, from (\ref{e30}) in Prop.~\ref{prop4}, $F_2(s)<F_1(s)$ when $s<s_n(\delta)$.\hfill\qedsymbol{}

\section{Illustration for $n=1$ and $n=3$} \label{sec7}
\mbox{} \\[-9mm]

To describe the illustration, we briefly summarize the results of Secs.~\ref{sec4}, \ref{sec5} and \ref{sec6}. For $n=0,1,..$ and $\delta\in(0,n+1)$, the unique positive solution of the equation
\beq \label{e178}
e^s=1+\frac{s^1}{1!}+...+\frac{s^n}{n!}+\frac{1}{\delta}~\frac{s^{n+1}}{n!}
\eq
is the unique maximizer over $s>0$ of
\beq \label{e179}
E_{n,\delta}(s)=\frac{1-\Bigl(1+\dfrac{s^1}{1!}+...+\dfrac{s^n}{n!}\Bigr)\,e^{-s}} {s^{\delta}}~.
\eq
For the maximum value $ME_{n,\delta}$ of $E_{n,\delta}(s)$, $s>0$, we have
\beq \label{e180}
ME_{n,\delta}=E_{n,\delta}(s_n(\delta))=F_1(s_n(\delta))=F_2(s_n(\delta))~,
\eq
where for $s>0$
\beq \label{e181}
F_1(s)=\frac{s^{n+1-\delta}}{\delta\cdot n!\,e^s}~,~~~~~~ F_2(s)= \frac{1}{\delta\cdot n!}~\frac{s^{n+1-\delta}} {1+\dfrac{s^1}{1!}+...+\dfrac{s^n}{n!}+\dfrac{1}{\delta}~\dfrac{s^{n+1}}{n!}}~.
\eq
When $S(\delta)$ is a lower bound for $s_n(\delta)$ with $S(\delta)>n+1-\delta$, we have
\beq \label{e182}
E_{n,\delta}(S(\delta))<ME_{n,\delta}<F_2(S(\delta))<F_1(S(\delta))~, ~~~~~~\delta\in(0,n+1)~.
\eq
For $n=0,1,...$ and $\delta\in(0,n+1)$, there is the bound
\beq \label{e183}
z+(n+1-\delta)<s_n(\delta)<z+(n+1-\delta)\,{\rm ln}\,z-(n+1)\, \frac{\partial\gamma}{\partial a}\,(1,z)~,
\eq
where $z={\rm ln}((n+1)/\delta)$ and
\beq \label{e184}
\frac{\partial\gamma}{\partial a}\,(1,z)=\il_0^z\,e^{-t}\,{\rm ln}\,t\,dt~, ~~~~~~z>0~.
\eq
Finally, the minimum value $\hat{E}_n$ of $ME_{n,\delta}$ as a function of $\delta\in(0,n+1)$ is assumed for
\beq \label{e185}
\delta=\delta_{n,E}=\Bigl(n!\Bigl(e-\sum_{k=0}^n\,\frac{1}{k!}\Bigr)\Bigr)^{-1}~,
\eq
for which we have $s_n(\delta_{n,E})=1$, and it equals
\beq \label{e186}
\hat{E}_n=\frac{1}{n!\,\delta_{n,E}\,e}=\frac1e\,\Bigl(e-\sum_{k=0}^n\,\frac{1}{k!} \Bigr)~.
\eq
There are the bounds
\begin{eqnarray} \label{e187}
& \mbox{} & \frac{n+1}{n+2}\,(n+1)<\delta_{n,E}<\frac{n+2}{n+3}\,(n+1)~, \nonumber \\[3mm]
& & \frac{n+3}{n+2}~\frac{1}{(n+1)!\,e}<\hat{E}_n<\frac{n+2}{n+1}~\frac{1}{(n+1)!\,e}~.
\end{eqnarray}

We shall illustrate all these results for the cases $n=1$ and $n=3$, where we take $\delta=\delta_{n,E}$, see (\ref{e185}). In particular, we shall verify numerically for these cases that Newton iteration for solving (\ref{e178}) gives the iterands
\begin{eqnarray} \label{e188}
s^{(0)} & = & {\rm ln}\Bigl(\frac{n+1}{\delta}\Bigr)+(n+1-\delta)~; \nonumber \\[3mm]
s^{(j+1)} & = & s-\left.\frac{e^s-1-\dfrac{s^1}{1!}-...-\dfrac{s^n}{n!}-\dfrac{1}{\delta}~\dfrac{s^{n+1}} {n!}}{e^s-1-\dfrac{s^1}{1!}-...-\dfrac{s^{n-1}}{(n-1)!}-\dfrac{n+1}{\delta} ~\dfrac{s^n}{n!}}\,\right|_{s=s^{(j)}}
\end{eqnarray}
with $j=0,1,...$ satisfying $s^{(j)}\pr1=s_n(\delta_{n,E})$ when $\delta=\delta_{n,E}$.

As to the upper bound in (\ref{e183}), we use the computational result (\ref{e158})
\beq \label{e189}
\frac{\partial\gamma}{\partial a}\,(1,z)=\Phi(z)=(1-e^{-z})\,{\rm ln}\,z+\sum_{k=1}^{\infty}\,\frac{({-}z)^k}{k\cdot k!}~,~~~~~~z>0~.
\eq
{\bf Case $n=1$}~~We have, see (\ref{e185}) and (\ref{e186}),
\beq \label{e190}
\delta=\delta_{1,E}=(e-2)^{-1}=1.392211191\,,~~\hat{E}_1=\frac1e\,(e-2)= 0.264241117~,
\eq
and, see (\ref{e187}),
\beq \label{e191}
\tfrac43<\delta<\tfrac32~,~~~~~~0.2452522960<\hat{E}_1<0.275909580~.
\eq
With $z={\rm ln}(2/\delta)=0.362253912$, we find for the lower bound in (\ref{e183}) and the value of $E_{1,\delta}$ at this lower bound
\beq \label{e192}
z+2-\delta=0.970042721=:S~,~~~~~~E_{1,\delta}(S)=0.264174903~.
\eq
Taking $s^{(0)}=0.970042721=S$, the Newton iteration
\beq \label{e193}
s^{(j+1)}=s-\frac{e^s-1-s-s^2/\delta}{e^s-1-2s/\delta}\Bigl|_{s=s^{(j)}}~,
\eq
gives for the iterands $s^{(0)}$, $s^{(1)}$, $s^{(2)}$, $s^{(3)}$ the respective values
\beq \label{e194}
0.970042721\,,~~1.002253487\,,~~1.000011471\,,~~1.000000000~.
\eq
We find for the upper bound in (\ref{e183})
\beq \label{e195}
z+(2-\delta)\,{\rm ln}\,z-2\,\frac{\partial\gamma}{\partial a}\,(1,z)= 1.026090795=: T
\eq
in which (\ref{e189}) has been used with 10 terms of the series for $5\times 10^{-10}$ absolute accuracy. We have
\beq \label{e196}
E_{1,\delta}(T)=0.264192948.
\eq
With $s=S$ equal to the lower bound in (\ref{e192}), we find
\beq \label{e197}
F_1(s)=0.267289754~,~~~~~~F_2(s)=0.26649408~.
\eq
\mbox{} \\
{\bf Case $n=3$}~~We have, see (\ref{e185}) and (\ref{e186}),
\beq \label{e198}
\delta=\delta_{3,E}=(6(e-\tfrac83))^{-1}=3.229025365~,~~~~~~\hat{E}_3= \frac1e\,(e-\tfrac83)=0.018988156~,
\eq
and, see (\ref{e187}),
\beq \label{e199}
\tfrac{16}{5}<\delta<\tfrac{10}{3}~,~~~~~~0.018393972<\hat{E}_3< 0.019160387~.
\eq
With $z={\rm ln}(4/\delta)=0.214114013$, we find for the lower bound in (\ref{e183}) and the value of $E_{3,\delta}$ at this lower bound
\beq \label{e200}
z+4-\delta=0.985088648=:S~,~~~~~~E_{3,\delta}(S)=0.018986579~.
\eq
Taking $s^{(0)}=0.985088648=S$, the Newton iteration
\beq \label{e201}
s^{(j+1)}=s- \frac{e^s-1-s-\tfrac12\,s^2-\tfrac16\,s^3-s^4/6\delta} {e^s-1-s-\tfrac12\,s^2-2s^3/3\delta}\Bigl|_{s=s^{(j)}}
\eq
gives for the iterands $s^{(0)}$, $s^{(1)}$, $s^{(2)}$, $s^{(3)}$ the respective values
\beq \label{e202}
0.985088648\,,~~1.001007241\,,~~1.000004222\,,~~1.000000000~.
\eq
We find for the upper bound in (\ref{e183})
\beq \label{203}
z+(4-\delta)\,{\rm ln}\,z-4\,\frac{\partial\gamma}{\partial a}\,(1,z)=1.026821936=: T
\eq
in which (\ref{e189}) has been used with 10 terms of the series for (more than) $5\times 10^{-10}$ absolute accuracy. We have
\beq \label{e204}
E_{3,\delta}(T)=0018983199~.
\eq
With $s=S$ equal to the lower bound in (\ref{e200}), we find
\beq \label{e205}
F_1(s)=0.019051464~,~~~~~~F_2(s)=0.019050286~.
\eq

\section{Basic properties of $G_{n,\delta}$ and $u_n$} \label{sec8}
\mbox{} \\[-9mm]

We recall that $G_{n,\delta}(u)$ is defined for $n=0,1,...\,$, $\delta\in[n,n+1]$ and $u>0$ by
\beq \label{e206}
G_{n,\delta}(u)=({-}1)^{n+1}\,\frac{e^{-u}-\Bigl(1-\dfrac{u^1}{1!}+...+({-}1)^n\, \dfrac{u^n}{n!}\Bigr)}{u^{\delta}}~,
\eq
and we define
\beq \label{e207}
G_{n,\delta}(0)=\dlim_{u\downarrow 0}\,G_{n,\delta}(u)=\left\{\ba{llll}
0 & \!\!, & ~~~\delta\in[n,n+1) & \!\!, \\[2mm]
\dfrac{1}{(n+1)!} & \!\!, & ~~~\delta=n+1 & \!\!.
\ea\right.
\eq
Then $G_{n,\delta}(u)$ is a continuous function of $u\geq0$. The case $n=0$ gives
\beq \label{e208}
G_{0,\delta}(u)=\frac{1-e^{-u}}{u^{\delta}}=E_{0,\delta}(u)~,~~~~~~ u\geq0~,
\eq
and has already been considered. We take $n=1,2,...$ in the sequel. From Taylor's theorem, we have for $u>0$ that there is a $\xi\in(0,u)$ such that
\beq \label{e209}
e^{-u}=1-\frac{u^1}{1}+...+({-}1)^n \,\frac{u^n}{n!}+({-}1)^{n+1}\,\frac{u^{n+1}}{(n+1)!}\,e^{-\xi}~.
\eq
Hence, $G_{n,\delta}(u)>0$ for $u>0$. We are interested in the maximum of $G_{n,\delta}(u)$ with $u\geq0$. When $\delta<n$, we have $\lim_{u\pr\infty}\,G_{n,\delta}(u)=\infty$, and when $\delta>n+1$, we have $\lim_{u\downarrow0}\,G_{n,\delta}(u)=\infty$. Therefore, we restrict ourselves to the case that $\delta\in[n,n+1]$. \\ \\
{\bf Proof of (\ref{e24}).}~~As to the first item in (\ref{e24}), we argue by induction. For $n=1$ and $u\geq0$, we have
\beq \label{e210}
e^{-u}-(1-u)\leq u
\eq
since $e^{-u}\leq1$, $u\geq0$, and
\beq \label{e211}
\lim_{u\pr\infty}\,\frac{e^{-u}-(1-u)}{u}=1~.
\eq
Hence,
\beq \label{e212}
0\leq G_{1,1}(u)\leq1\,,~~u\geq0~;~~~~~~\lim_{u\pr\infty}\,G_{1,1}(u)=1~,
\eq
and this is the case $n=1$ in the first item of (\ref{e24}). Assume now that the first item in (\ref{e24}) holds for some $n=1,2,...\,$. We want to show that
\beq \label{e213}
({-}1)^{n+2}\Bigl(e^{-u}-\Bigl(1-\frac{u^1}{1!}+...+({-}1)^{n+1}\, \frac{u^{n+1}} {(n+1)!}\Bigr)\Bigr)\leq \frac{u^{n+1}}{(n+1)!}~,~~~~~~u\geq0~.
\eq
For $u=0$, there is equality in (\ref{e213}), and so it is sufficient to show that
\begin{eqnarray} \label{e214}
& \mbox{} & \frac{d}{du}\,\Bigl[({-}1)^{n+2}\Bigl(e^{-u}-\Bigl( 1-\frac{u^1}{1!}+ \frac{u^2}{2!}-...+({-}1)^{n+1}\,\frac{u^{n+1}}{(n+1)!}\Bigr)\Bigr)\Bigr] \nonumber \\[3mm]
& & \leq~\frac{d}{du}\,\Bigl[\frac{u^{n+1}}{(n+1)!}\Bigr]~,~~~~~~ u\geq0~.
\end{eqnarray}
The right-hand side of (\ref{e214}) equals $u^n/n!$ while the left-hand side of (\ref{e214}) equals
\begin{eqnarray} \label{e215}
& \mbox{} & ({-}1)^{n+2}\Bigl({-}e^{-u}-\Bigl({-}1+u-...+({-}1)^{n+1}\,\frac{u^n}{n!} \Bigr)\Bigr) \nonumber \\[3mm]
& & =~({-}1)^{n+1}\Bigl(e^{-u}-\Bigl(1-\frac{u^1}{1!}+...+({-}1)^n\,\frac{u^n}{n!} \Bigr)\Bigr)\leq\frac{u^n}{n!}~,
\end{eqnarray}
where the inequality holds by the induction hypothesis. Obviously, we also have
\beq \label{e216}
\lim_{u\pr\infty}\,\frac{({-}1)^{n+2}\Bigl(e^{-u}-\Bigl(1-\dfrac{u^1}{1!}+...+({-}1) ^{n+1}\,\dfrac{u^{n+1}}{(n+1)!}\Bigr)\Bigr)}{u^{n+1}}=\frac{1}{(n+1)!}~,
\eq
and so the first item in (\ref{e24}) holds for $n+1$.

As to the second item in (\ref{e24}), we use Taylor's theorem, see (\ref{e209}), and we see that the inequality in the second item of (\ref{e24}) holds since $\xi>0$ in (\ref{e209}). Moreover, since $\xi\pr0$ as $u\downarrow 0$, the limit relation in the second item of (\ref{e24}) holds as well.\hfill\qedsymbol \\
\mbox{}

\noindent The observation that $G_{n,n}(u)=\frac{1}{n!}-G_{n-1,n}(u)$ and $G_{n-1,n}(u)>0$ for $u>0$ gives an alternative proof of the first item in (\ref{e24}). \\\mbox{}

We restrict now to the cases that $\delta\in(n,n+1)$ for which we have
\beq \label{e217}
\lim_{u\pr\infty}\,G_{n,\delta}(u)=0=\lim_{u\downarrow0}\,G_{n,\delta}(u)~.
\eq
{\bf Proof of Proposition~\ref{prop13}.}~~Let $n=1,2,...$ and $\delta\in(n,n+1)$. Then for $u>0$, we have $G_{n,\delta}'(u)=0$ if and only if
\begin{eqnarray} \label{e218}
& \mbox{} & \Bigl({-}e^{-u}-\Bigl({-}1+\frac{u^1}{1!}-...-({-}1)^{n-1}\, \frac{u^{n-1}}{(n-1)!}\Bigr)\Bigr)\,u \nonumber \\[3mm]
& & -~\delta\Bigl(e^{-u}-\Bigl(1-\frac{u^1}{1!}+\frac{u^2}{2!}-...+({-}1)^n\,\frac{u^n}{n!} \Bigr)\Bigr)=0~,
\end{eqnarray}
i.e., if and only if
\begin{eqnarray} \label{e219}
-e^{-u}(u+\delta) & \!+ & \!1\cdot(u+\delta)-\frac{u^1}{1!}\,(u+\delta)+...+({-}1)^{n-1} \,\frac{u^{n-1}}{(n-1)~}\,(u+\delta) \nonumber \\[3mm]
& \!+ & \!\delta({-}1)^n\,\frac{u^n}{n!}=0~.
\end{eqnarray}
This gives (\ref{e52}).

We next show that the equation (\ref{e53}) has exactly one positive solution $u>0$, and to that end, we consider the function
\beq \label{e220}
K_{n,\delta}(u)=({-}1)^n\,\Bigl[e^{-u}-\Bigl(1-\frac{u^1}{1!}+...+({-}1)^{n-1}\, \frac{u^{n-1}}{(n-1)!}+\frac{({-}1)^n\,\delta\,u^n}{n!(u+\delta)}\Bigr)\Bigr]
\eq
for $u\geq0$. For small $u>0$, we have
\begin{eqnarray} \label{e221}
& \mbox{} & \hspace*{-4mm}K_{n,\delta}(u) \nonumber \\[3mm]
& & \hspace*{-4mm}=~({-}1)^n\,\Bigl[\frac{({-}1)^n}{n!}\,u^n+ \frac{({-}1)^{n+1}}{(n+1)!}\,u^{n+1}+\frac{({-}1)^{n+2}}{(n+2)!}\,u^{n+2}+...- \frac{({-}1)^n\,u^n}{n!(1+u/\delta)}\Bigr] \nonumber \\[3mm]
& & \hspace*{-4mm}=~\frac{u^n}{n!}-\frac{u^{n+1}}{(n+1)!}+\frac{u^{n+2}}{(n+2)!}-...- \frac{u^n}{n!}\,\Bigl(1-\frac{u}{\delta}+\frac{u^2}{\delta^2}-...\Bigr) \nonumber \\[3mm]
& & \hspace*{-4mm}=~\Bigl(\frac{1}{\delta}-\frac{1}{n+1}\Bigr)\,\frac{u^{n+1}}{n!}- \Bigl(\frac{1}{\delta^2}-\frac{1}{(n+1)(n+2)}\Bigr)\,\frac{u^{n+2}}{n!}+...~.
\end{eqnarray}
Writing $K=K_{n,\delta}$ for brevity, we thus have
\beq \label{e222}
K(0)=K'(0)=...=K^{(n)}(0)=0~;~~~~~~K^{(n+1)}(0)=\frac{n+1}{\delta}-1>0
\eq
since $\delta\in(n,n+1)$. Furthermore, for large $u>0$
\begin{eqnarray} \label{e223}
K(u) & = & {-}\,\frac{\delta u^n}{n!(u+\delta)}+\frac{u^{n-1}}{(n-1)!}+ O(u^{n-2}) \nonumber \\[3mm]
& = & {-}\,\frac{u^{n-1}}{n!}\,\Bigl(\frac{\delta u}{u+\delta}-n\Bigr) +O(u^{n-2})={-}\,\frac{u^{n-1}}{n!}\,(\delta-n)+O(u^{n-2})~. \nonumber \\
\mbox{}
\end{eqnarray}
Hence, $K(u)<0$ for large $u>0$ since $\delta\in(n,n+1)$.

We next write
\beq \label{e224}
\frac{u^n}{u+\delta}=\frac{(u+\delta-\delta)^n}{u+\delta}= \frac{({-}\delta)^n}{u+\delta}+\sum_{k=1}^n\,\Bigl(\!\ba{c} n \\[1mm] k\ea\!\Bigr)(u+\delta)^{k-1}({-}\delta)^{n-k}~,
\eq
and we compute
\begin{eqnarray} \label{e225}
K^{(n)}(u) & = & ({-}1)^n\,\Bigl[({-}1)^n\,e^{-u}- \frac{({-}1)^n\,\delta}{n!} \cdot({-}\delta)^n\cdot\frac{{-}1\cdot{-}2\cdot...\cdot{-}n}{(u+\delta)^{n+1}} \Bigr] \nonumber \\[3mm]
& = & e^{-u}-\Bigl(\frac{1}{1+u/\delta}\Bigr)^{n+1}~.
\end{eqnarray}
The function
\beq \label{e226}
u\geq0\mapsto \varp(u):=u-(n+1)\,{\rm ln}\Bigl(1+\frac{u}{\delta}\Bigr)
\eq
is strictly convex, and we have $\varp(0)=0$, $\varp'(0)=1-(n+1)/\delta<0$ since $\delta\in(n,n+1)$, and $\varp(u)\pr\infty$ as $u\pr\infty$. Therefore, $\varp(u)$ has a unique positive zero. We have
\beq \label{e227}
\varp(u)={\rm ln}\,[e^u/(1+u/\delta)^{n+1}]=0\Leftrightarrow K^{(n)}(u)=0
\eq
so that $\varp$ and $K^{(n)}$ have the same positive zeros. Denote the unique positive zero of $\varp$ and $K^{(n)}$ by $k_n$.

From $K^{(n-1)}(0)=0$, we infer that $K^{(n-1)}(u)$ is positive and strictly increasing in $u\in(0,k_n)$, and strictly decreasing in $u\in(k_n,\infty)$. We have $\lim_{u\pr\infty}\,K^{(n-1)}(u)<0$, for otherwise $K^{(n-1)}(u)\geq0$ for all $u>0$ and so $K^{(n-2)}(u)>0,K^{(n-3)}(u)>0,...,K^{(0)}(u)=K(u)>0$ for all $u>0$, contradicting that $K(u)<0$ for large $u$. Therefore, $K^{(n-1)}$ has a unique positive zero $k_{n-1}>k_n$, with $K^{(n-1)}(u)>0$ when $u\in(0,k_{n-1})$ and $K^{(n-1)}(u)<0$ when $u\in(k_{n-1},\infty)$. Continuing this way, we find that $K^{(n)},K^{(n-1)},...,K^{(1)},K^{(0)}$ have unique positive zeros $k_n,k_{n-1},...,k_1,k_0$ with $0<k_n<k_{n-1}<...<k_1<k_0$ and $K^{(j)}(u)>0$ when $0<u<k_j$ while $K^{(j)}(u)<0$ when $u>k_j$, $j=n,n-1,...,0\,$. This proves the existence and uniqueness of a positive zero of $K=K_{n,\delta}$.\hfill\qedsymbol{} \\ \\
{\bf Consequence of the proof of Prop.~\ref{prop13}.}~~We have that $K_{n,\delta}(u)$ is strictly decreasing in $\delta\in(n,n+1)$ for any $u>0$. We also have that $\frac{d}{du}\,[K_{n,\delta}(u)]<0$ at $u=k_0$ (since $k_1<k_0$). Therefore, $k_n$ is strictly decreasing in $\delta\in(n,n+1)$. \\ \\
{\bf Proof of Proposition~\ref{prop14}.}~~Let $n=1,2,...$ and $\delta\in(n,n+1)$, and denote the unique positive solution $u$ of (\ref{e53}) by $u_n(\delta)$. We consider
\beq \label{e228}
L(\delta,u)=K_{n,\delta}(u)~,~~~~~~u\geq0~,
\eq
so that $L(\delta,u_n(\delta))=0$ for $\delta\in(n,n+1)$. Writing $u=u_n(\delta)$ for brevity, we have
\beq \label{e229}
\frac{\partial L}{\partial\delta}\,(\delta,u)={-}\,\frac{u^{n+1}}{n!}~ \frac{1}{(u+\delta)^2}~,
\eq
\begin{eqnarray} \label{e230}
& \mbox{} & \hspace*{-7mm}\frac{\partial L}{\partial u}\,(\delta,u_n(\delta)) \nonumber \\[3mm]
& & \hspace*{-7mm}=~{-}({-}1)^n\Bigl[e^{-u}-\Bigl(1-\frac{u^1}{1!}+...+ ({-}1)^{n-2}\,\frac{u^{n-2}}{(n-2)!}\Bigr) \nonumber \\[3mm]
& & \hspace*{2.2cm}-~({-}1)^n\, \frac{\delta}{n!}~\frac{(n-1)\,u^n+n\delta u^{n-1}}{(u+\delta)^2}\Bigr] \nonumber \\[3mm]
& & \hspace*{-7mm}=~{-}({-}1)^n\Bigl[({-}1)^{n-1}\, \frac{u^{n-1}}{(n-1)!}+ ({-}1)^n\,\frac{\delta u^n}{n!\,(u+\delta)}\Bigr]-\frac{\delta}{n!}~ \frac{(n-1)\,u^n+n\delta u^{n-1}}{(u+\delta)^2} \nonumber \\[3mm]
& & \hspace*{-7mm}=~\frac{u^{n-1}}{(n-1)!}-\frac{\delta u^n}{n!(u+\delta)}- \frac{\delta}{n!}~\frac{(n-1)\,u^n+n\delta u^{n-1}}{(u+\delta)^2}~,
\end{eqnarray}
where we have used that $u=u_n(\delta)$ satisfies (\ref{e53}). Then from $0=L(\delta,u_n(\delta))$, we get
\begin{eqnarray} \label{e231}
0 & \!\!= & \!\!\frac{d}{d\delta}\,[L(\delta,u_n(\delta))]=\frac{\partial L}{\partial\delta}\,(\delta,u_n(\delta))+\frac{\partial L}{\partial u}\,(\delta,u_n(\delta))\,u_n'(\delta) \nonumber \\[3mm]
& \!\!= & \!\!{-}\,\frac{u^{n+1}}{n!}~\frac{1}{(u+\delta)^2}+u' \Bigl(\frac{u^{n-1}}{(n-1)!}-\frac{\delta\,u^n}{n!(u+\delta)}- \frac{\delta}{n!}~ \frac{(n-1)\,u^n+n\delta u^{n-1}}{(u+\delta)^2}\Bigr) \nonumber \\[3mm]
& \!\!= & \!\!\frac{u^{n+1}}{n!}\,\Bigl({-}\,\frac{u^2}{(u+\delta)^2}+u'\, \frac{(n-\delta)\,u^2+(n+1-\delta)\,\delta u}{(u+\delta)^2}\Bigr)
\end{eqnarray}
in which we have written $u=u_n(\delta)$ in the last two lines of (\ref{e231}). Hence,
\beq \label{e232}
u=u'((n-\delta)\,u+(n+1-\delta)\,\delta)~.
\eq
It has been observed after the proof of Prop.~\ref{prop13} that $u'=u_n'(\delta)<0$. Since $u=u_n(\delta)>0$, it thus follows that $(\delta-n)\,u-(n+1-\delta)\,\delta>0$ and this yields (\ref{e55}) in Prop.~\ref{prop14}.\hfill\qedsymbol{} \\ \\
{\bf Proof of Proposition~\ref{prop15}.}~~The quantity in the second member of the two equivalence propositions in (\ref{e56}) and (\ref{e57}) coincides with $K_{n,\delta}(u)$ in (\ref{e220}). From the proof of Prop.~\ref{prop13} we have for $u>0$ that $u<u_n(\delta)\Leftrightarrow K_{n,\delta}(u)>0$ and $u>u_n(\delta)\Leftrightarrow K_{n,\delta}(u)<0$.\hfill\qedsymbol{} \\
\mbox{}

We now address the problem of computing $u_n(\delta)$ for $n=1,2,...$ and $\delta\in(n,n+1)$. From the proof of Prop.~\ref{prop13} we have
\beq \label{e233}
K_{n,\delta}(u)=0~,~~~~~~u=u_n(\delta)=k_0
\eq
and
\beq \label{e234}
K_{n,\delta}'(u)<0\,,~~u>k_1~;~~~~~~K_{n,\delta}''(u)<0\,,~~u>k_2~,
\eq
where $k_2<k_1<k_0$. To find $u_n(\delta)$, we must solve equation (\ref{e53}) for $u$, or equivalently, the equation $K_{n,\delta}(u)=0$ with $K_{n,\delta}$ given in (\ref{e220}). We can use for this Newton's method, where we choose $u^{(0)}$ equal to an upper bound for $u_n(\delta)$, causing the Newton iterands (on account of (\ref{e234})) to converge monotonically to $u_n(\delta)$. The Newton iteration step is given by
\beq \label{e235}
u^{(j+1)}=u^{(j)}-\frac{K_{n,\delta}(u^{(j)})}{K_{n,\delta}'(u^{(j)})}~, ~~~~~~j=0,1,...~,
\eq
in which
\begin{eqnarray} \label{e236}
& \mbox{} & \hspace*{-7mm}-\,\frac{K_{n,\delta}(u)}{K_{n,\delta}'(u)} \nonumber \\[2mm]
& & \hspace*{-7mm}=~\frac{e^{-u}-1+\dfrac{u^1}{1!}-...-({-}1)^{n-1}\,\dfrac{u^{n-1}} {(n-1)!}-({-}1)^n\,\dfrac{\delta u^n}{n!(\delta+u)}}
{e^{-u}-1+\dfrac{u^1}{1!}-...-({-}1)^{n-2}\,\dfrac{u^{n-2}}{(n-2)!}-({-}1)^{n-1}\, \dfrac{\delta u^{n-1}((n-1)u+n\delta)}{n!(\delta+u)^2}}~. \nonumber \\
\mbox{}
\end{eqnarray}
For $\delta$ not too close to $n+1$ and $n$ not too large there are no real problems. However, when $y=n+1-\delta$ is small, we should operate more carefully. For small values of $y$, one can compute $u_n(\delta)$, alternatively, by using the B\"urmann-Lagrange inversion formula. Thus, one writes the equation (\ref{e53}) as
\beq \label{e237}
\frac{\delta}{u+\delta}=\frac{({-}1)^n\,n!}{u^n}\,\Bigl(e^{-u}-1+\frac{u^1}{1!} -...-({-}1)^{n-1}\,\frac{u^{n-1}}{(n-1)!}\Bigr)=:g(u)~.
\eq
With $\delta=n+1-y$ and solving for $y$, one gets
\beq \label{e238}
\frac{u\,g(u)+(n+1)(g(u)-1)}{g(u)-1}=y~.
\eq
From $\exp({-}u)=\sum_{k=0}^{\infty}\,({-}1)^k\,u^k/k!$, it is found that
\beq \label{e239}
g(u)=\frac{({-}1)^n\,n!}{u^n}\,\sum_{k=n}^{\infty}\,({-}1)^k\,\frac{u^k}{k!} =\sum_{k=0}^{\infty}\,({-}1)^k\,\frac{n!}{(n+k)!}\,u^k~,
\eq
\beq \label{e240}
u\,g(u)+(n+1)(g(u)-1)=\sum_{k=1}^{\infty}\,({-}1)^k\, \frac{k\cdot n!} {(n+k+1)!}\,u^{k+1}~,
\eq
\beq \label{e241}
g(u)-1=\sum_{k=1}^{\infty}\,({-}1)^k\,\frac{n!}{(n+k)!}\,u^k~.
\eq
Some further manipulations then bring (\ref{e238}) in the form
\beq \label{e242}
u~\frac{\dsum_{k=0}^{\infty}\,e_k\,u^k} {\dsum_{k=0}^{\infty}\,d_k\,u^k} =(n+2)\,y~,
\eq
where
\beq \label{e243}
e_k=({-}1)^k\,\frac{(k+1)(n+2)!}{(n+k+2)!}~,~~~~~~d_k=({-}1)^k\, \frac{(n+1)!} {(n+k+1)!}~,~~~~~~k=0,1,...~.
\eq
Thus (compare formula (\ref{e104}) for the $S$-case), one has for small $y>0$
\beq \label{e244}
u=\sum_{m=1}^{\infty}\,h_m((n+2)\,y)^m~,
\eq
where
\beq \label{e245}
h_m=\frac1m\,C_{u^{m-1}}\left[\left( \frac{\dsum_{k=0}^{\infty}\,d_k\,u^k} {\dsum_{k=0}^{\infty}\,e_k\,u^k}\right)^m\right]~,~~~~~~m=1,2,...~.
\eq
Then we find, for example, with $y=2-\delta$,
\beq \label{e246}
u_1(y)=3y+\tfrac16\,(3y)^2+\tfrac{2}{45}\,(3y)^3+\tfrac{7}{540}\,(3y)^4+\tfrac {113}{28350}\,(3y)^5+...
\eq
and, with $y=4-\delta$,
\beq \label{e247}
u_3(y)=5y+\tfrac{2}{15}\,(5y)^2+\tfrac{38}{1575}\,(5y)^3+\tfrac{439}{94500} \,(5y)^4+\tfrac{9131}{9922500}\,(5y)^5+...~.
\eq
\mbox{}

We owe the reader still a proof of the last identity in (\ref{e54}).

\begin{lem} \label{lem5}
Let $n=1,2,...$ and $\delta\in(n,n+1)$. Then we have
\beq \label{e248}
MG_{n,\delta}=\max_{u\geq0}\,G_{n,\delta}(u)=G_{n,\delta}(u_n(\delta)) =H_{n,\delta}(u_n(\delta))~,
\eq
where
\beq \label{e249}
H_{n,\delta}(u)=\frac{1}{n!}~\frac{u^{n+1-\delta}}{u+\delta}~,~~~~~~ u\geq0~.
\eq
\end{lem}

\noindent
{\bf Proof.}~~The first identity in (\ref{e248}) is just the definition of $MG_{n,\delta}$ in (\ref{e54}). Since $u_n(\delta)$ is the unique maximizer of $G_{n,\delta}(u)$ over $u\geq0$, the second identity in (\ref{e248}) is also obvious. The third identity in (\ref{e248}) follows from the definition of $G_{n,\delta}(u)$ in (\ref{e20}) and the fact that $u_n(\delta)$ is the unique positive solution of (\ref{e53}), together with a simple computation.\hfill\qedsymbol{}

\section{Bounds for $u_n(\delta)$} \label{sec9}
\mbox{} \\[-9mm]

According to Taylor's theorem with remainder in integral form, see (\ref{e111}), applied to $f(u)=\exp({-}u)$ and with $n-1$ instead of $n$, we have for $u>0$
\beq \label{e250}
({-}1)^n\Bigl(e^{-u}-\Bigl( 1-\frac{u^1}{1!}+...+({-}1)^{n-1}\,\frac{u^{n-1}}{(n-1)!}\Bigr)\Bigr)= \frac{1}{(n-1)!}\,e^{-u}\,\il_0^u\,x^{n-1}\,e^x\,dx~.
\eq
When we combine this with Prop.~\ref{prop15}, we see that for $n=1,2,...\,$, $\delta\in(n,n+1)$ and $u>0$
\begin{eqnarray} \label{e251}
u\:<\:{\rm or}\:>\:u_n(\delta) & \Leftrightarrow & \frac{1}{(n-1)!}\,e^{-u}\,\il_0^u\, x^{n-1}\,e^x\,dx\:>\:{\rm or}\:<\: \frac{\delta u^n}{n!(u+\delta)} \nonumber \\[3mm]
& \Leftrightarrow & \il_0^u\,x^{n-1}\,e^x\,dx\:>\:{\rm or}\:<\:\frac1n~\frac{\delta u^n}{u+\delta}\,e^u~.
\end{eqnarray}
{\bf Proof of Proposition\ref{prop16}.} \mbox{}

a.~~Assume that $U\in C^1(n,n+1]$ with $U(n+1)=0<U(\delta)$, $\delta\in(n,n+1)$. According to (\ref{e251}), we have for $\delta\in(n,n+1)$ that $U(\delta)<u_n(\delta)$ if and only if
\beq \label{e252}
\il_0^{U(\delta)}\,x^{n-1}\,e^x\,dx>\frac1n~\frac{\delta U^n(\delta)} {U(\delta)+\delta}\,e^{U(\delta)}~.
\eq
We have equality in (\ref{e252}) when $\delta=n+1$ since $U(n+1)=0$. Hence, it is sufficient to show that for $\delta\in(n,n+1)$
\beq \label{e253}
\frac{d}{d\delta}\,\left[\il_0^{U(\delta)}\,x^{n-1}\,e^x\,dx\right]< \frac{d}{d\delta}\,\Bigl[\frac1n~\frac{\delta U^n(\delta)}{U(\delta)+\delta}\,e^{U(\delta)}\Bigr]~.
\eq
We compute for $\delta\in(n,n+1)$
\beq \label{e254}
\frac{d}{d\delta}\,\left[\il_0^{U(\delta)}\,x^{n-1}\,e^x\,dx\right]= U'(\delta)\,U^{n-1}(\delta)\,e^{U(\delta)}~.
\eq
Furthermore, for $\delta\in(n,n+1)$,
\begin{eqnarray} \label{e255}
& \mbox{} & \frac{d}{d\delta}\,\Bigl[\frac{\delta U^n(\delta)}{U(\delta)+\delta}\,e^{U(\delta)}\Bigr] = \frac{\delta U^n(\delta)}{U(\delta)+\delta}\,U'(\delta)\,e^{U(\delta)} \nonumber \\[3mm]
& & \hspace*{1.5cm}+~\frac{(U^n(\delta)+n\delta U^{n-1}(\delta)\,U'(\delta))(U(\delta)+\delta)-\delta U^n(\delta)(U'(\delta)+1)} {(U(\delta)+\delta)^2} \nonumber \\[3mm]
& & =~U^{n-1}\,e^U\,\Bigl[\frac{\delta UU'}{U+\delta}+ \frac{(U+n\delta U')(U+\delta)-\delta U(U'+1)} {(U+\delta)^2}\Bigr] \nonumber \\[3mm]
& & =~U^{n-1}\,e^U\,\Bigl[\frac{\delta UU'}{U+\delta}+\frac{n\delta(U+\delta)-\delta U}{(U+\delta)^2}\,U'+\frac{U^2}{(U+\delta)^2}\Bigr]~,
\end{eqnarray}
where we have written in the last two lines $U=U_n(\delta)$ and $U'=U_n'(\delta)$ for brevity. Hence, (\ref{e253}) holds for $\delta\in(n,n+1)$ when
\beq \label{e256}
nU'<\frac{\delta UU'}{U+\delta}+\frac{n\delta(U+\delta)-\delta U}{(U+\delta)^2}\, U'+\frac{U^2}{(U+\delta)^2}~,
\eq
i.e., when
\beq \label{e257}
U'\,\frac{n(U+\delta)^2-\delta U(U+\delta)-(n-1)\,\delta U-n\delta^2} {(U+\delta)^2}<\frac{U^2}{(U+\delta)^2}~,
\eq
i.e., when
\beq \label{e258}
U'({-}(\delta-n)\,U^2+\delta U(n+1-\delta))<U^2~,
\eq
i.e., when
\beq \label{e259}
U(\delta)>{-}U'(\delta)((\delta-n)\,U(\delta)-(n+1-\delta)\,\delta)~,
\eq
where we have used that $U(\delta)>0$. This is (\ref{e61}).

b.~~The proof of b is completely similar to the proof of a when $U(n+1)=0$, with reversed inequality signs. In the case that $U(n+1)>0$, we need to show that
\beq \label{e260}
\il_0^{U(\delta)}\,x^{n-1}\,e^x\,dx\leq\frac{\delta U^n(\delta)} {n(U(\delta)+\delta)}\,e^{U(\delta)}
\eq
holds for $\delta=n+1$, and then we can proceed as in the proof of a, with reversed inequality signs. This point is settled by Lemma~\ref{lem6} below.\hfill\qedsymbol{}

\begin{lem} \label{lem6}
We have for $n>0$ and $U\geq0$
\beq \label{e261}
\il_0^U\,x^{n-1}\,e^x\,dx\leq\frac{(n+1)\,U^n}{n(U+n+1)}\,e^U~,
\eq
with equality if and only if $U=0$.
\end{lem}

\noindent
{\bf Proof.}~~We have equality in (\ref{e261}) when $U=0$. Hence, it is sufficient to show that for $U>0$
\beq \label{e262}
\frac{d}{dU}\,\left[\il_0^U\,x^{n-1}\,e^x\,dx\right]<\frac{n+1}{n}~\frac{d}{dU}
\,\Bigl[\frac{U^n\,e^U} {U+n+1}\Bigr]~,
\eq
i.e., that
\beq \label{e263}
U^{n-1}\,e^U<\frac{n+1}{n}\,U^{n-1}\,e^U\Bigl(\frac{U}{U+n+1}+\frac{(n-1)\,U+n(n+1)} {(U+n+1)^2}\Bigr)~,
\eq
i.e., that
\beq \label{e264}
1<\frac{n+1}{n}~\frac{U^2+2nU+n(n+1)}{(U+n+1)^2}~.
\eq
The inequality in (\ref{e264}) is equivalent with
\beq \label{e265}
\frac{n}{n+1}\,U^2+2nU+n(n+1)<U^2+2nU+n(n+1)~,
\eq
and this obviously holds for $n>0$ and $U>0$.\hfill\qedsymbol \\
\mbox{}

We shall now present the proofs of Props.~\ref{prop17}, \ref{prop18} and \ref{prop19} as consequences of Prop.~\ref{prop16}. \\ \\
{\bf Proof of Proposition~\ref{prop17}.}~~Let $n=1,2,...\,$, and set
\beq \label{e266}
U(\delta)=\frac{(n+1-\delta)\,\delta}{\delta-n}+(n+1-\delta)~,~~~~~~ \delta\in(n,n+1]~.
\eq
We have $U(n+1)=0$, and we shall check that
\beq \label{e267}
U(\delta)>{-}U'(\delta)((\delta-n)\,U(\delta)-(n+1-\delta)\,\delta)~, ~~~~~~\delta\in(n,n+1)~,
\eq
i.e., that
\beq \label{e268}
\frac{\delta}{\delta-n}+1>{-}U'(\delta) \Bigl((\delta-n)\Bigl(\frac{\delta} {\delta-n}+1\Bigr)-\delta\Bigr)={-}U'(\delta)(\delta-n)~,~~~~~~\delta\in(n,n+1)~.
\eq
Now $U'(\delta)={-}n/(\delta-n)^2-2$, so checking (\ref{e268}) amounts to checking that $2\delta-n>n+2(\delta-n)^2$, i.e., that $\delta-n>(\delta-n)^2$. The latter inequality indeed holds since $\delta\in(n,n+1)$. Hence, the first inequality in (\ref{e63}) of Prop.~\ref{prop17} follows from Prop.~\ref{prop16}a.

Next, let $n=2,3,...\,$, and set
\beq \label{e269}
U(\delta)=\frac{\delta}{\delta-n}~,~~~~~~\delta\in(n,n+1]~.
\eq
We have $U(n+1)=n+1\geq0$, and we shall check that
\beq \label{e270}
U(\delta)<{-}U'(\delta)((\delta-n)\,U(\delta)-(n+1-\delta)\,\delta)~,~~~~~~ \delta\in(n,n+1)~,
\eq
i.e., that
\begin{eqnarray} \label{e271}
& \mbox{} & \frac{\delta}{\delta-n}<{-}U'(\delta)\Bigl((\delta-n)\,\frac{\delta}{\delta-n}-(n+1-\delta) \,\delta\Bigr)={-}U'(\delta)\,\delta(\delta-n)~, \nonumber \\[2mm]
& & \hspace*{8.5cm}\delta\in(n,n+1)~.
\end{eqnarray}
Since $U'(\delta)={-}n/(\delta-n)^2$, checking (\ref{e271}) amounts to checking that $1<n$. Hence, (\ref{e271}) holds since $n=2,3,...\,$, and the second inequality in (\ref{e63}) follows for this case from Prop.~\ref{prop16}b.

For $n=1$, we have that (\ref{e271}) does not hold, and to handle this case, we use criterion (\ref{e57}) in Prop.~\ref{prop15}. Thus, for $u>0$
\beq \label{e272}
u>u_1(\delta)\Leftrightarrow e^{-u}>1-\frac{\delta u}{u+\delta}
\eq
when $\delta\in(1,2)$. Now $1-\delta u/(u+\delta)=0$ for $u=\delta/(\delta-1)$, and so $\delta/(\delta-1)>u_1(\delta)$ when $\delta\in(1,2)$.\hfill\qedsymbol{} \\
\mbox{}

As a step towards the sharp results in Props.~\ref{prop18} and \ref{prop19} we first prove the inequality (\ref{e67}) that is useful due to its simplicity and sharpness, especially when $\delta\uparrow n+1$.

\begin{lem} \label{lem7}
Let $n=1,2,...\,$. Then
\beq \label{e273}
u_n(\delta)<\frac{(n+1-\delta)\,\delta}{\delta-n}-{\rm ln}(\delta-n)~,~~~~~~ \delta\in(n,n+1)~.
\eq
\end{lem}

\noindent
{\bf Proof.}~~Set
\beq \label{e274}
U(\delta)=\frac{(n+1-\delta)\,\delta}{\delta-n}-{\rm ln}(\delta-n)~,~~~~~~ \delta\in(n,n+1]~.
\eq
We have $U(n+1)=0$, and we shall show that
\beq \label{e275}
U(\delta)<U'(\delta)(\delta-n)\,{\rm ln}(\delta-n)~,~~~~~~ \delta\in(n,n+1)~,
\eq
which is what the condition in Prop.~\ref{prop16}b amounts to for this $U$. We have $U'(\delta)={-}n/(\delta-n)^2-1-1/(\delta-n)$, and so we must verify that for $\delta\in(n,n+1)$
\beq \label{e276}
\frac{(n+1-\delta)\,\delta}{\delta-n}-{\rm ln}(\delta-n)<{-}\Bigl(\frac{n}{\delta-n}+(\delta-n)+1\Bigr)\,{\rm ln}(\delta-n)~,
\eq
i.e., that
\beq \label{e277}
((\delta-n)^2+n)\cdot{-}{\rm ln}(\delta-n)>(n+1-\delta)\,\delta~.
\eq
Now, for $\delta\in(n,n+1)$,
\beq \label{e278}
-{\rm ln}(\delta-n)={-}{\rm ln}(1-(n+1-\delta))>(n+1-\delta)(1+\tfrac12\,(n+1-\delta))~.
\eq
Hence, it suffices to show that for $\delta\in(n,n+1)$
\beq \label{e279}
((\delta-n)^2+n)(1+\tfrac12\,(n+1-\delta))>\delta~,
\eq
i.e., that
\beq \label{e280}
\tfrac12\,(n+1-\delta)((\delta-n)^2+n)>\delta-n-(\delta-n)^2 =(\delta-n)(n+1-\delta)~,
\eq
i.e., that
\beq \label{e281}
\tfrac12\,(\delta-n)^2-(\delta-n)+\tfrac12\,n=\tfrac12\,(\delta-n-1)^2+ \tfrac12\,(n-1)>0~.
\eq
Since $n\geq1$ and $\delta<n+1$, the latter inequality indeed holds, and so Prop.~\ref{prop16}b gives the result.\hfill\qedsymbol{} \\ \\
{\bf Proof of Proposition~\ref{prop18}.}~~Let $n=2,3,...\,$, and set
\beq \label{e282}
U(\delta)=\frac{(n+1-\delta)\,\delta}{\delta-n}+(n+1-\delta)+\tfrac12\,(n+1-\delta)^2 ~,~~~~~~\delta\in(n,n+1)~.
\eq
We have $U(n+1)=0$, and we shall verify condition (\ref{e62}) in Prop.~\ref{prop16}b. We compute $U'(\delta)={-}n/(\delta-n)^2-2-(n+1-\delta)$, and so we must check for $\delta\in(n,n+1)$ that
\begin{eqnarray} \label{e283}
& \mbox{} & \frac{(n+1-\delta)\,\delta}{\delta-n}+(n+1-\delta)+\tfrac12\,(n+1-\delta)^2 \nonumber \\[3mm]
& & <~\Bigl(\frac{n}{(\delta-n)^2}+2+(n+1-\delta)\Bigr)(\delta-n)((n+1-\delta)+\tfrac12 \,(n+1-\delta)^2)~, \nonumber \\
\mbox{}
\end{eqnarray}
i.e., that
\beq \label{e284}
\frac{\delta}{\delta-n}+1+\tfrac12\,(n+1-\delta)<\Bigl(\frac{n}{(\delta-n)^2}+2 +(n+1-\delta)\Bigr)(\delta-n)(1+\tfrac12\,(n+1-\delta))~,
\eq
i.e., that
\beq \label{e285}
2\delta-n+\tfrac12\,(\delta-n)(n+1-\delta)< (n+(2+(n+1-\delta))(\delta-n)^2) (1+\tfrac12\,(n+1-\delta))~.
\eq
Set $y=n+1-\delta\in(0,1)$, so that $\delta=n+1-y$ and $\delta-n=1-y$. We should verify that for $y\in(0,1)$
\beq \label{e286}
n+2-2y+\tfrac12\,(1-y)\,y<(n+(2+y)(1-y)^2)(1+\tfrac12\,y)~,
\eq
i.e., that
\beq \label{e287}
-4y+(1-y)\,y<ny+(2+y)^2(1-y)^2-4=ny-(y+y^2)(4-y-y^2)~,
\eq
i.e., that
\beq \label{e288}
-3-y<n-4-3y+2y^2+y^3~,
\eq
i.e., that
\beq \label{e289}
1<n-2y+2y^2+y^3~.
\eq
We have $n=2,3,...$ and $y\in(0,1)$ and so $n-2y+2y^2\geq 3/2$. Therefore, (\ref{e289}) indeed holds, and we get from Prop.~\ref{prop16}b
\beq \label{e290}
u_n(\delta)<\frac{(n+1-\delta)\,\delta}{\delta-n}+(n+1-\delta)+ \tfrac12\,(n+1-\delta)^2~,~~~~~~\delta\in(n,n+1)~.
\eq

Next, still with $n=2,3,...\,$, set
\beq \label{e291}
U(\delta)=\frac{(n+1-\delta)\,\delta}{\delta-n}+1~,~~~~~~\delta\in(n,n+1]~.
\eq
We have $U(n+1)=1\geq0$, and we shall verify condition (\ref{e62}) in Prop.~\ref{prop16}b, i.e., that for $\delta\in(n,n+1)$
\beq \label{e292}
U(\delta)<{-}U'(\delta)(\delta-n)~.
\eq
We have $U'(\delta)={-}n/(\delta-n)^2-1$, and so we should verify that for $\delta\in(n,n+1)$
\beq \label{e293}
(n+1-\delta)\,\delta+(\delta-n)<n+(\delta-n)^2~,
\eq
i.e., that
\beq \label{e294}
2\delta^2-(3n+2)\,\delta+n^2+2n>0~.
\eq
The left-hand side of (\ref{e294}) equals $(\delta-n)(2\delta-(n+2))>0$ when $\delta\in(n,n+1)$ and $n=2,3,...\,$. We therefore get from Prop.~\ref{prop16}b
\beq \label{e295}
u_n(\delta)<\frac{(n+1-\delta)\,\delta}{\delta-n}+1~,~~~~~~ \delta\in(n,n+1)~.
\eq
Together with (\ref{e290}) we then get (\ref{e65}).\hfill\qedsymbol{} \\ \\
It is observed that both (\ref{e289}) and (\ref{e294}) do not hold for all $y=n+1-\delta$ with $\delta\in(n,n+1)$ when $n=1$. We have that $u_1(1.25)$ exceeds both the right-hand side of (\ref{e290}) and the right-hand side of (\ref{e295}) when $n=1$ and $\delta=1.25$. We do have (\ref{e290}) when $n=1$ and $\delta\downarrow1$ and we do have (\ref{e295}) when $n=1$ and $\delta\uparrow2$. \\ \\
{\bf Outline of the proof of Proposition~\ref{prop19}.}~~With
\beq \label{e296}
U(\delta)=\frac{(2-\delta)\,\delta}{\delta-1}+(2-\delta)+\tfrac12\,(2-\delta)^2 +\tfrac13\,(2-\delta)^3~,~~~~~~\delta\in(1,2]~,
\eq
we have $U(2)=0$, and we should verify the condition (\ref{e62}) in Prop.~\ref{prop16}b with $n=1$ and this $U$. This goes in the usual way. With $y=2-\delta\in(0,1)$, we have to check at an intermediate stage validity of
\beq \label{e297}
2(1-y)-\tfrac12\,y^2-\tfrac13\,y^3<(2+y+y^2)(1+\tfrac12\,y+\tfrac13\,y^2)(1-y^2)~,
\eq
followed by a check at the final stage of validity of
\beq \label{e298}
-\,\tfrac23+\tfrac76\,y-\tfrac56\,y^2<\tfrac16\,y^3+\tfrac13\,y^4~.
\eq
The left-hand side of (\ref{e298}) is less than $-\frac46+\frac86\,y-\frac46\,y^2$ which equals $-\,\frac23\,(1-y)^2<0$ for $y\in(0,1)$.\hfill\qedsymbol \\
\mbox{}

With $y=n+1-\delta$, there are the following expansions
\begin{eqnarray} \label{e299}
u_n(y) & \!\!= & \!\!(n+2)\,y+\frac{(n+1)(n+2)}{n+3}\,y^2 \nonumber \\[3mm]
& & \!\!+\,\Bigl(\frac{2n+2}{(n+2)(n+3)^2(n+4)}\,+\,\Bigl(\frac{n+1}{(n+2)(n+3)} \Bigr)^2\Bigr)\,(n+2)^2\,y^2{+}...\,, \nonumber \\
\mbox{}
\end{eqnarray}
\beq \label{e300}
\frac{(n+1-\delta)\,\delta}{\delta-n}=(n+1)\,y+ny^2+ny^3+...~,
\eq
\beq \label{e301}
\frac{(n+1-\delta)\,\delta}{\delta-n}-{\rm ln}(\delta-n)=(n+2)\,y+(n+\tfrac12) \,y^2+(n+\tfrac13)\,y^3+...~,
\eq
\beq \label{e302}
\frac{\delta}{\delta-n}=n+1+ny+ny^2+ny^3+...~,
\eq
\beq \label{e303}
\frac{(n+1-\delta)\,\delta}{\delta-n}+(n+1-\delta)=(n+2)\,y+ny^2+ny^3+...~,
\eq
\beq \label{e304}
\frac{(n+1-\delta)\,\delta}{\delta-n}+(n+1-\delta)+\tfrac12\,(n+1-\delta)^2=(n+2)\,y+(n+\tfrac12)\, y^2+ny^3+...~.
\eq
For the case $n=1$, we have with $y=2-\delta$
\beq \label{e305}
\frac{(2-\delta)\,\delta}{\delta-1}+(2-\delta)+\tfrac12\,(2-\delta)^2 +\tfrac13\,(2-\delta)^3=3y+\tfrac32\,y^2+\tfrac43\,y^3+y^4+y^5+...\,,
\eq
\beq \label{e306}
u_1(\delta)=3y+\tfrac32\,y^2+\tfrac65\,y^3+\tfrac{21}{20}\,y^4+\tfrac{339}{350}\,y^5+...~.
\eq
The expansions in (\ref{e299}) and (\ref{e306}) follow from the B\"urmann-Lagrange theorem, see (\ref{e242})--(\ref{e247}) in Sec.~\ref{sec8}. \\ \\
{\bf Proof of Proposition~\ref{prop20}.}~~Let $\delta\in(1,2)$. We have that $u_1(\delta)$ is the solution $u>0$ of the equation
\beq \label{e307}
e^{-u}=1-\frac{\delta u}{u+\delta}~,
\eq
which can be rewritten as
\beq \label{e308}
u=B\Bigl(1-\frac{u}{e^u-1}\Bigr)~,~~~~~~B=\frac{\delta}{\delta-1}~.
\eq
We have by Prop.~\ref{prop17}
\beq \label{309}
B>u_1(\delta)>\frac{(2-\delta)\,\delta}{\delta-1}+(2-\delta)=B-2(\delta-1)~.
\eq
Hence, with $u=u_1(\delta)$,
\beq \label{e310}
\frac{B}{e^B-1}<\frac{u}{e^u-1}<\frac{B-2(\delta-1)}{e^{B-2(\delta-1)}-1}
\eq
since $u/(e^u-1)$ is a decreasing function of $u\geq0$. Both the first and the third member in (\ref{e310}) are $O(B\,e^{-B})$ as $\delta\downarrow1$, and this gives (\ref{e69}).\hfill\qedsymbol{} \\
\mbox{}

We conclude this section by showing that $u_n(\delta)$ is a strictly convex function of $\delta\in(n,n+1)$ where $n=1,2,...\,$. The corresponding result for $s_n$ required a particular lower bound for $s_n(\delta)$, $\delta\in(0,n+1)$, see Lemma~\ref{lem1} in Sec.~\ref{sec4} and the proof of Prop.~\ref{prop7}. The inequalities that we already have for the $u_n(\delta)$, $\delta\in(n,n+1)$, are so sharp that no additional effort is required. \\ \\
{\bf Proof of Proposition~\ref{prop21}.}~~Let $n=1,2,...\,$. We shall show that $u_n''(\delta)>0$, $\delta\in(n,n+1)$. Writing $u=u_n(\delta)$, $u'=u_n'(\delta)$ with $\delta\in(n,n+1)$ for brevity, we have by Prop.~\ref{prop14} and \ref{prop17}
\beq \label{e311}
u'={-}\,\frac{u}{(\delta-n)\,u-(n+1-\delta)\,\delta}~,
\eq
\beq \label{e312}
u>\frac{(n+1-\delta)\,\delta}{\delta-n}+(n+1-\delta)=(n+1-\delta)\Bigl( 1+\frac{\delta}{\delta-n}\Bigr)~.
\eq
Since the function $-x/((\delta-n)\,x-(n+1-\delta)\,\delta)$ is strictly decreasing in $x\in((n+1-\delta)\,\delta/(\delta-n),\infty)$, we have
\beq \label{e313}
u'>{-}\,\frac{(n+1-\delta)(1+\delta/(\delta-n))} {(n+1-\delta)\,\delta+(n+1-\delta)(\delta-n)-(n+1-\delta)\,\delta}={-}\, \frac{2\delta-n}{(\delta-n)^2}~.
\eq
We compute for $\delta\in(n,n+1)$ from (\ref{e311})
\beq \label{e314}
u''(\delta)={-}\,\frac{u'((\delta-n)\,u-(n+1-\delta)\,\delta)-u(u+(\delta-n)\,u'-(n+1)+2\delta)} {((\delta-n)\,u-(n+1-\delta)\,\delta)^2}~.
\eq
Hence, $u''(\delta)>0$ if and only if
\beq \label{e315}
C(\delta):=u^2+(n+1-\delta)\,\delta u'+(2\delta-n-1)\,u>0~.
\eq
For $\delta\in(n,n+1)$ and $n=1,2,...$
\beq \label{e316}
2\delta-n-1>2n-(n+1)=n-1\geq0~,
\eq
and so, using (\ref{e312}) and (\ref{e313}),
\begin{eqnarray} \label{e317}
& \mbox{} & \hspace*{-7mm}C(\delta) \nonumber \\[3mm]
& & \hspace*{-7mm}>~(n{+}1{-}\delta)^2\Bigl(1+\frac{\delta}{\delta-n}\Bigr)^2 -(n{+}1{-}\delta)\,\delta\,\frac{2\delta-n}{(\delta-n)^2} \nonumber \\[3mm]
& & \hspace*{4.2cm}+~(2\delta{-}n{-}1)(n{+}1{-}\delta)\Bigl(1+\frac{\delta} {\delta-n}\Bigr) \nonumber \\[3mm]
& & \hspace*{-7mm}=~\frac{n{+}1{-}\delta}{(\delta-n)^2}\,\left((n{+}1{-}\delta)
(2\delta-n)^2-\delta (2\delta-n)+(2\delta{-}n{-}1)(2\delta-n)(\delta-n)\right) \nonumber \\[3mm]
& & \hspace*{-7mm}=~\frac{(n{+}1{-}\delta)(2\delta-n)}{(\delta-n)^2}\,\left((n{+}1{-}
\delta)(2\delta-n)-\delta+ (2\delta{-}n{-}1)(\delta-n)\right)=0~. \nonumber \\
\mbox{}
\end{eqnarray}
Hence, $u_n''(\delta)>0$ for $\delta\in(n,n+1)$.\hfill\qedsymbol{}

\section{The function $MG_{n,\delta}$} \label{sec10}
\mbox{} \\[-9mm]

We consider in this section for $n=1,2,...$ the function
\beq \label{e318}
MG_{n,\delta}=G_{n,\delta}(u_n(\delta))=\frac{1}{n!}~\frac{u^{n+1-\delta}}{u+\delta} ~,~~~~~~u=u_n(\delta)~,
\eq
with $\delta\in[n,n+1]$. \\ \\
{\bf Proof of Proposition~\ref{prop22}.}~~We have for $\delta\in(n,n+1)$
\beq \label{e319}
\frac{\delta}{\delta-n}-\Bigl(\frac{(n+1-\delta)\,\delta}{\delta-n}+(n+1-\delta)\Bigr) =2\delta-(n-1)\in(n-1,n+1)~,
\eq
and so, by Prop.~\ref{prop17},
\beq \label{e320}
u_n(\delta)=\frac{\delta}{\delta-n}\,(1+O(\delta-n))~,~~~~~~\delta\downarrow n~.
\eq
Therefore, from (\ref{e318}),
\begin{eqnarray} \label{e321}
MG_{n,\delta} & = & \frac{1}{n!}~\frac{\Bigl(\dfrac{\delta}{\delta-n}\Bigr)^{n+1-\delta}} {\dfrac{\delta}{\delta-n}+\delta}\,(1+O(\delta-n)) \nonumber \\[3mm]
& = & \frac{1}{n!}~\frac{\delta^{n-\delta}(\delta-n)^{\delta-n}} {1+(\delta-n)}\,(1+O(\delta-n)) \nonumber \\[3mm]
& = & \frac{1}{n!}\,\Bigl(1+O\Bigl((\delta-n)\,{\rm ln}\Bigl(\frac{1}{\delta-n} \Bigr)\Bigr)\Bigr)~,~~~~~~\delta\downarrow n~.
\end{eqnarray}
Next, from the first inequality in Prop.~\ref{prop17} and Lemma~\ref{lem7}
\beq \label{e322}
u_n(\delta)=(n+2)(n+1-\delta)(1+O(n+1-\delta))~,~~~~~~\delta\uparrow n+1~,
\eq
where we have used that $-{\rm ln}(\delta-n)=(n+1-\delta)+O((n+1-\delta)^2)$. Therefore, from (\ref{e318}),
\begin{eqnarray} \label{e323}
MG_{n,\delta} & = & \frac{1}{n!}~\frac{(n+2)^{n+1-\delta}} {(n+1)+O(n+1-\delta)}\,(n+1-\delta)^{n+1-\delta}(1+O(n+1-\delta)) \nonumber \\[3mm]
& = & \frac{1}{(n+1)!}\,\Bigl(1+O\Bigl((n+1-\delta)\,{\rm ln}\Bigl( \frac{1}{n+1-\delta}\Bigr)\Bigr)\Bigr)~,~~~~~~\delta\uparrow n+1~, \nonumber \\
\mbox{}
\end{eqnarray}
as required.\hfill\qedsymbol{} \\ \\
We conclude that $MG_{n,\delta}$ is a continuous function of $\delta\in[n,n+1]$. \\ \\
{\bf Proof of Proposition~\ref{prop23}.}~~Writing $u=u_n(\delta)$, $u'=u_n'(\delta)$, we have from Prop.~\ref{prop14} and (\ref{e318})
\begin{eqnarray} \label{e324}
& \mbox{} & \hspace*{-7mm}\frac{d}{d\delta}\,[{\rm ln}\,MG_{n,\delta}] \nonumber \\[3mm]
& & \hspace*{-7mm}=~\frac{d}{d\delta}\,[(n+1-\delta)\,{\rm ln}\,u_n(\delta)-{\rm ln}(u_n(\delta)+\delta)] \nonumber \\[3mm]
& & \hspace*{-7mm}=~{-}{\rm ln}\,u+(n+1-\delta)\,\frac{u'}{u}-\frac{u'+1}{u+\delta} \nonumber \\[3mm]
& & \hspace*{-7mm}=~{-}{\rm ln}\,u-\frac{n+1-\delta}{(\delta-n)\,u-(n{+}1{-}\delta)\,\delta} +\frac{u}{(u+\delta)((\delta-n)\,u-(n{+}1{-}\delta)\,\delta)}-\frac{1}{u+\delta} \nonumber \\[3mm]
& & \hspace*{-7mm}=~{-}{\rm ln}\,u+\frac{u-(n+1-\delta)(u+\delta)} {(u+\delta)((\delta-n)\,u-(n+1-\delta)\,\delta)}-\frac{1}{u+\delta} \nonumber \\[3mm]
& & \hspace*{-7mm}=~{-}{\rm ln}\,u+\frac{(\delta-n)\,u-(n+1-\delta)\,\delta} {(u+\delta)((\delta-n)\,u-(n+1-\delta)\,\delta)}-\frac{1}{u+\delta}={-}{\rm ln}\,u~,
\end{eqnarray}
as required.\hfill\qedsymbol{} \\ \\
{\bf Proof of Proposition~\ref{prop24}.}~~The function $u_n(\delta)$ is continuous in $\delta\in(n,n+1]$ and decreases strictly from $+\infty$ at $\delta=n$ to 0 at $\delta=n+1$. Hence, there is a unique $\delta=\delta_{n,G}$ such that $u_n(\delta)=1$. Thus, we have
\begin{eqnarray} \label{e325}
\hat{G}_n & = & \min_{\delta\in(n,n+1)}\,MG_{n,\delta}= MG_{n,\delta_{n,G}}=G_{n,\delta_{n,G}}(u_n(\delta_{n,G})) \nonumber \\[3mm]
& = & G_{n,\delta_{n,G}}(1)=({-}1)^{n+1}\,\left. \frac{e^{-u}-\Bigl(1-\dfrac{u^1}{1!}+...+({-}1)^n\,\dfrac{u^n}{n!}\Bigr)} {u^{\delta_{n,G}}}\right|_{u=1} \nonumber \\[3mm]
& = & ({-}1)^{n+1}\Bigl(e^{-1}-\Bigl(1-\frac{1}{1!}+...+\frac{({-}1)^n} {n!}\Bigr)\Bigr)~,
\end{eqnarray}
and this gives (\ref{e75}). Next, from (\ref{e318}) with $u=u_n(\delta_{n,G})=1$,
\beq \label{e326}
\hat{G}_n=\frac{1}{n!}~\frac{1}{1+\delta_{n,G}}~,~~~~~~ \delta_{n,G}=(n!\,\hat{G}_n)^{-1}-1~,
\eq
and this gives (\ref{e76}).\hfill\qedsymbol{} \\
\mbox{}

We finally give some bounds for $\delta_{n,G}$ and $\hat{G}_n$.

\begin{lem} \label{lem8}
Let $n=1,2,...\,$. Then
\beq \label{e327}
n+1-\frac{1}{n+2}<\delta_{n,G}<n+1~,~~~~~~\frac{n+1}{n+2}~\frac{1}{(n+1)!} <\hat{G}_n<\frac{n+2}{n+3}~\frac{1}{(n+1)!}~.
\eq
\end{lem}

\noindent
{\bf Proof.}~~We have from (\ref{e75})
\begin{eqnarray} \label{e328}
n!\,\hat{G}_n & = & ({-}1)^{n+1}\,n!\Bigl(e^{-1}-\sum_{k=0}^n\,\frac {({-}1)^k}{k!}\Bigr) \nonumber \\[3mm]
& = & ({-}1)^{n+1}\,\sum_{k=n+1}^{\infty}\,\frac {({-}1)^k\,n!}{k!} \nonumber \\[3mm]
& = & \frac{1}{n+1}-\frac{1}{(n+1)(n+2)}+\frac{1}{(n+1)(n+2)(n+3)}-...~.
\end{eqnarray}
Therefore,
\beq \label{e329}
n!\,\hat{G}_n>\frac{1}{n+1}-\frac{1}{(n+1)(n+2)}=\frac{1}{n+2}~,
\eq
\begin{eqnarray} \label{e330}
n!\,\hat{G}_n & < & \frac{1}{n+1}-\frac{1}{(n+1)(n+2)}+\frac{1} {(n+1)(n+2)(n+3)} \nonumber \\[3mm]
& = & \frac{1}{n+2}+\frac{1}{(n+1)(n+2)(n+3)} \nonumber \\[3mm]
& = & \frac{1}{n+2}~\frac{(n+1)(n+3)+1}{(n+1)(n+3)}=\frac{n+2}{(n+1)(n+3)}~,
\end{eqnarray}
and this gives (\ref{e327}).\hfill\qedsymbol{}

\newpage
\section{Lower and upper bounds for $MG_{n,\delta}$ and illustration of the results for the cases $n=1$ and $n=3$} \label{sec11}
\mbox{} \\[-9mm]

Let $n=1,2,...$ and let $U(\delta)$ be an approximation of $u_n(\delta)$, $\delta\in(n,n+1)$. Then $G_{n,\delta}(U(\delta))$ is a lower bound for $MG_{n,\delta}$. For the case that $U(\delta)$ is a lower bound for $u_n(\delta)$ such that $U(\delta)\geq(n+1-\delta)\,\delta/(\delta-n)$, $\delta\in(n,n+1)$, it has been argued at the end of Subsec.~\ref{subsec2.2} that
\beq \label{e331}
MG_{n,\delta}<H(U(\delta))~,~~~~~~\delta\in(n,n+1)~,
\eq
where for $\delta\in(n,n+1)$
\beq \label{e332}
H(u)=\frac{1}{n!}~\frac{u^{n+1-\delta}}{u+\delta}~,~~~~~~ u\geq0~.
\eq
This result depends on the following property of $H$.

\begin{lem} \label{lem9}
Let $\delta\in(n,n+1)$. The function $H(u)$ is strictly increasing in $u\in(0,(n+1-\delta)\,\delta/(\delta-n))$ and strictly decreasing in $u\in((n+1-\delta)\,\delta/(\delta-n),\infty)$.
\end{lem}

\noindent{\bf Proof.}~~We compute for $u>0$
\beq \label{e333}
\frac{d}{du}\,\Bigl[\frac{u^{n+1-\delta}}{u+\delta}\Bigr]=u^{n-\delta}\, \frac{(n+1-\delta)\,\delta-(\delta-n)\,u}{(u+\delta)^2}~,
\eq
and this is positive for $u\in(0,(n+1-\delta)\,\delta/(\delta-n))$ and negative for $u\in((n+1-\delta)\,\delta/(\delta-n),\infty)$.\hfill\qedsymbol{} \\
\mbox{}

We next illustrate the results of Secs.~\ref{sec9} and \ref{sec10} and those on the lower and upper bounds for $MG_{n,\delta}$ for the cases that $n=1$ and $n=3$. We first summarize the main results. For $n=1,2,...$ and $\delta\in(n,n+1)$, the unique positive solution $u=u_n(\delta)$ of the equation
\beq \label{e334}
e^{-u}=1-\frac{u^1}{1!}+...+({-}1)^{n-1}\,\frac{u^{n-1}}{(n-1)!}+({-}1)^n \,\frac{\delta u^n}{n!(\delta+u)}
\eq
is the unique positive maximizer of the quantity
\beq \label{e335}
G_{n,\delta}(u)=({-}1)^{n+1}\,\frac{e^{-u}-1+\dfrac{u^1}{1!}-\dfrac{u^2}{2!}+...-({-}1)^n \,\dfrac{u^n}{n!}}{u^{\delta}}~,~~~~~~u>0~.
\eq
For the maximal value $MG_{n,\delta}$ of $G_{n,\delta}(u)$ we have
\beq \label{e336}
MG_{n,\delta}=G_{n,\delta}(u_n(\delta))=H(u_n(\delta))~,
\eq
where $H$ is given by (\ref{e332}). The minimal value $\hat{G}_n$ of $MG_{n,\delta}$ as a function of $\delta\in(n,n+1)$ is assumed for
\beq \label{e337}
\delta=\delta_{n,G}=\Bigl[({-}1)^{n+1}\,n!\Bigl(e^{-1}-\Bigl( 1-\frac{1}{1!}+...+({-}1)^n\,\frac{1}{n!}\Bigr)\Bigr)\Bigr]^{-1}-1~,
\eq
that satisfies $u_n(\delta_{n,G})=1$, and it equals
\beq \label{e338}
\hat{G}_n=\frac{1}{n!(1+\delta_{n,G})}=({-}1)^{n+1}\Bigl( e^{-1}-\Bigl(1-\frac{1}{1!}+...+({-}1)^n\,\frac{1}{n!}\Bigr)\Bigr)~.
\eq
There are the bounds
\beq \label{e339}
n+1-\frac{1}{n+2}<\delta_{n,G}<n+1~,~~~~~~\frac{n+1}{n+2}~\frac{1}{(n+1)!} <\hat{G}_n<\frac{n+2}{n+3}~\frac{1}{(n+1)!}~.
\eq
For $u_1$ we have the bounds, see Props.~\ref{prop17} and \ref{prop19},
\beq \label{e340}
u_1(\delta)>\frac{(2-\delta)\,\delta}{\delta-1}+(2-\delta)~,~~~~~~ u_1(\delta)<\frac{\delta}{\delta-1}~,~~~~~~\delta\in(1,2)~,
\eq
and
\beq \label{e341}
u_1(\delta)<\frac{(2-\delta)\,\delta}{\delta-1}+(2-\delta)+\tfrac12\,(2-\delta)^2+\tfrac13\, (2-\delta)^3~,~~~~~~\delta\in(1,2)~.
\eq
For $u_3$, we have the bounds, see Props.~\ref{prop17} and \ref{prop18},
\beq \label{e342}
u_3(\delta)>\frac{(4-\delta)\,\delta}{\delta-3}+(4-\delta)~,~~~~~~ u_3(\delta)<\frac{\delta}{\delta-3}~,~~~~~~\delta\in(3,4)~,
\eq
and
\beq \label{e343}
u_3(\delta)<\frac{(4-\delta)\,\delta}{\delta-3}+\min\,\{1,4-\delta+\tfrac12\,(4-\delta)^2 \}~,~~~~~~\delta\in(3,4)~.
\eq

We shall take for both cases $\delta=\delta_{n,G}$, see (\ref{e337}). In particular, we shall verify numerically for these cases that the Newton method to solve the equation (\ref{e334}) has the iterands
\begin{eqnarray} \label{e344}
& \mbox{} & \hspace*{-7mm}u^{(0)}=\frac{(n+1-\delta)\,\delta}{\delta-n}+(n+1-\delta)~, \nonumber \\[3mm]
& & \hspace*{-7mm}u^{(j+1)} \nonumber \\[1mm]
& & \hspace*{-7mm}=~u+\left.\frac{e^{-u}-1+\dfrac{u^1}{1!}-...- ({-}1)^{n-1}\, \dfrac{u^{n-1}}{(n-1)!}-({-}1)^n\,\dfrac{\delta u^n}{n!(\delta+u)}} {e^{-u}{-}1{+}\dfrac{u^1}{1}{-}...{-}({-}1)^{n-2}\,\dfrac{u^{n-2}}{(n{-}2)!}{-}
({-}1)^{n-1}\, \dfrac{\delta u^{n-1}((n{-}1)u{+}n\delta)}{n!(\delta+u)^2}}\right|_{u=u^{(j)}} \nonumber \\[2mm]
\mbox{}
\end{eqnarray}
with $j=0,1,...$ that converge to $u_n(\delta_{n,G})=1$ when $\delta=\delta_{n,G}$. \\ \\
{\bf Case $n=1$.}~~We have, see (\ref{e337}) and (\ref{e338}),
\beq \label{e345}
\delta=\delta_{1,G}=[1\cdot1(e-(1-1))]^{-1}-1=e-1=1.718281828~,
\eq
\beq \label{e346}
\hat{G}_1=\frac{1}{1\cdot(1+e-1)}=\frac1e=0.367879441~,
\eq
and, see (\ref{e339}),
\beq \label{e347}
\tfrac53<\delta_{1,G}<2~,~~~~~~\tfrac13<\hat{G}_1<\tfrac38~.
\eq
We find for the lower bound in (\ref{e340}) and the values of $G_{1,\delta}$ at this lower bound, see (\ref{e335}),
\beq \label{e348}
\frac{(2-\delta)\,\delta}{\delta-1}+(2-\delta)=0.95564753=:U~,~~~~~~ G_{1,\delta}(U)=0.367791633~.
\eq
Taking $u^{(0)}=U=0.955647534$, the Newton iteration
\beq \label{e349}
u^{(j+1)}=u+\frac{e^{-u}-1+\delta u/(\delta+u)}{e^{-u}-(\delta/(\delta+u))^2}\Bigl|_{u=u^{(j)}}~,~~~~~~ j=0,1,...~,
\eq
gives for the iterands $u^{(0)}$, $u^{(1)}$, $u^{(2)}$, $u^{(3)}$ the respective values
\beq \label{e350}
0.955647534\,,~~1.002603361\,,~~1.000078440\,,~~1.000000000~.
\eq
We find for the upper bounds in (\ref{e340}) and (\ref{e341})
\beq \label{e351}
\dfrac{\delta}{\delta-1}=2.392211191~,
\eq
\beq \label{e352}
\frac{(2-\delta)\,\delta}{\delta-1}+(2-\delta)+\tfrac12\, (2-\delta)^2+\tfrac13\,(2-\delta)^3=1.002782965=: V~,
\eq
\beq \label{e353}
G_{1,\delta}(V)=0.367879108~.
\eq
With $u=U$ equal to the lower bound for $u_1(\delta)$ in (\ref{e348}), we find
\beq \label{e354}
H(u)=\frac{u^{2-\delta}}{u+\delta}=0.369232215~,
\eq
see (\ref{e331}--\ref{e332}). \\ \\
{\bf Case $n=3$.}~~We have, see (\ref{e337}) and (\ref{e338})
\beq \label{e355}
\delta=\delta_{3,G}=[6(e^{-1}-(1-1+\tfrac12-\tfrac16))]^{-1}{-}1=[6(e^{-1}-\tfrac13 )]^{-1}{-}1=3.824470167~,
\eq
\beq \label{e356}
\hat{G}_3=\frac{1}{6(1+\delta)}=0.034546107~,
\eq
and, see (\ref{e339}),
\beq \label{e357}
\tfrac{19}{5}<\delta<4~,~~~~~~\tfrac{1}{30}<\hat{G}_3<\tfrac{5}{144}~.
\eq
We find for the lower bound in (\ref{e342}) and the value of $G_{3,\delta}$ at this lower bound, see (\ref{e335}),
\beq \label{e358}
\frac{(4-\delta)\,\delta}{\delta-3}+(4-\delta)=0.989760152=:U~,~~~~~~ G_{3,\delta}(U)=0.034545828~.
\eq
Taking $u^{(0)}=U=0.989760152$, the Newton iteration
\beq \label{e359}
u^{(j+1)}=u+\frac{e^{-u}-1+u-\frac12\,u^2+\frac16\,\delta u^3/(u+\delta)} {e^{-u}-1+u-\frac16\,\delta u^2(2u+3\delta)/(u+\delta)^2}\Bigl|_{u=u^{(j)}}~,~~~~~j=0,1,...\,,
\eq
gives for the iterands $u^{(0)}$, $u^{(1)}$, $u^{(2)}$, $u^{(3)}$ the respective values
\beq \label{e360}
0.989760152\,,~~1.000383438\,,~~1.000000511\,,~~1.000000000~.
\eq
We find for the upper bounds in (\ref{e342}) and (\ref{e343})
\beq \label{e361}
\frac{\delta}{\delta-3}=4.638700487~,
\eq
\beq \label{e362}
\frac{(4-\delta)\,\delta}{\delta-3}+(4-\delta)+\tfrac12\,(4-\delta)^2= 1.005165513=: V~,
\eq
\beq \label{e363}
G_{3,\delta}(V)=0.034546037~.
\eq
With $u=U$ equal to the lower bound for $u_3(\delta)$ in (\ref{e358}), we find
\beq \label{e364}
H(u)=\frac{1}{3!}~\frac{u^{4-\delta}}{u+\delta}=0.034557097~,
\eq
see (\ref{e331}--\ref{e332}).

\end{document}